\documentclass{elsarticle}
\usepackage{amsfonts, cleveref}
\usepackage{amsmath,amssymb, amsthm}
\usepackage[all]{xy}
\usepackage{tikz}
\usepackage{enumitem,units,multirow}
\usepackage{boldline,colortbl}
\usepackage{mathtools}
\usepackage{hhline,url}
\usetikzlibrary{trees}
\usetikzlibrary{arrows,shapes,snakes,automata,backgrounds,petri}
\usepackage{accents}
\usepackage{bm}
\newlength{\dhatheight}
\newcommand{\doublehat}[1]{%
    \settoheight{\dhatheight}{\ensuremath{\hat{#1}}}%
    \addtolength{\dhatheight}{-0.35ex}%
    \hat{\vphantom{\rule{1pt}{\dhatheight}}%
    \smash{\hat{#1}}}}

\mathchardef\mhyphen="2D 

\makeatletter
\newcommand{\specialnumber}[1]{%
  \def\tagform@##1{\maketag@@@{(\ignorespaces##1\unskip\@@italiccorr#1)}}%
}
\newcommand{\specialeqref}[2]{\begingroup
  \def\tagform@##1{\maketag@@@{(\ignorespaces##1\unskip\@@italiccorr#2)}}%
  \eqref{#1}\endgroup}
\makeatother

\newtheorem{thm}{Theorem}
\newtheorem{thm1}{Theorem}[section]
\newtheorem{lemma1}[thm1]{Lemma}
\newtheorem{conj}[thm1]{Conjecture}

\newtheorem{def1}[thm1]{Definition}
\newtheorem{example}[thm1]{Example}
\newtheorem{notation}[thm1]{Notation}
\newtheorem{prop}[thm1]{Proposition}

\crefformat{appendix}{#2#1#3}

\theoremstyle{remark}
\newtheorem{remark}[thm1]{Remark}

\newproof{pf}{Proof}
\newenvironment{myproof}[1] {{\em Proof of {#1}. }}{\hfill$\square$}

\definecolor{light-gray}{gray}{0.9}
\definecolor{mygray}{gray}{0.8}

\begin{document}
\begin{frontmatter}

\title{Inheritance of oscillation in chemical reaction networks}

\author[ref1]{Murad Banaji\corref{cor1}}\address[ref1]{Middlesex University, London, Department of Design Engineering and mathematics, The Burroughs, London NW4 4BT, UK.} \cortext[cor1]{{\tt m.banaji@mdx.ac.uk}}

\begin{abstract}
Some results are presented on how oscillation is inherited by chemical reaction networks (CRNs) when they are built in natural ways from smaller oscillatory networks. The main results describe four important ways in which a CRN can be enlarged while preserving its capacity for oscillation. The results are for general CRNs, not necessarily fully open, but lead to an important corollary for fully open networks: if a fully open CRN $\mathcal{R}$ with mass action kinetics admits a nondegenerate (resp., linearly stable) periodic orbit, then so do all such CRNs which include $\mathcal{R}$ as an induced subnetwork. This claim holds for other classes of kinetics, but fails, in general, for CRNs which are not fully open. Where analogous results for multistationarity can be proved using the implicit function theorem alone, the results here call on regular and singular perturbation theory. Equipped with these results and with the help of some analysis and numerical simulation, lower bounds are put on the proportion of small fully open CRNs capable of stable oscillation under various assumptions on the kinetics. This exploration suggests that small oscillatory motifs are an important source of oscillation in CRNs. 
\end{abstract}

\end{frontmatter}

\section{Introduction and context of the paper}

At the heart of many biological systems are chemical reaction networks (CRNs), and the question of when these admit oscillation is of both theoretical and practical interest. Oscillation is known to occur -- and play a key role -- in a great variety of biological contexts. Examples include the natural rhythms of body clocks and ovulation, biochemical oscillations in cellular signalling, cyclic behaviour of various diseases, and periodic fluctuations in Lotka-Volterra-type models of interacting populations. Several chapters of \cite{MurrayMathBio} and \cite{mathphys} detail mathematical models of oscillation in biological settings. Some general biological principles underlying biological oscillation are discussed in \cite{novaktyson}. Once a network admitting oscillation is identified, we might naturally wonder whether this network occurs as a ``motif'' in other larger networks and, if so, whether the larger networks must themselves admit oscillation. The desire to phrase this question precisely and provide some simple and partial answers motivates this work.

Several papers have treated analogous questions about the inheritance of multistationarity in CRNs \cite{joshishiu,feliuwiufInterface2013,Joshi.2013aa,JoshiShiu2016}. In a recent contribution it was shown that a great deal can be done in this direction using the implicit function theorem \cite{banajipanteaMPNE}. An (incomplete) list of network modifications proven to preserve the property of admitting nondegenerate multistationarity were listed; these collectively define a partial order $\preceq$ on the set of all CRNs such that if a CRN $\mathcal{R}$ admits nondegenerate multistationarity, then so do all CRNs $\succeq \mathcal{R}$ in this partial order. Although it is likely that most, if not all, of the results in \cite{banajipanteaMPNE} can be restated with ``nondegenerate oscillation'' replacing ``nondegenerate multistationarity'', only part of this task is undertaken here: we prove four results about general CRNs, Theorems~\ref{thmnewdepreac}~to~\ref{thmnewwithopen}, which are analogues of related results about multistationarity in \cite{banajipanteaMPNE}, also numbered Theorems~1~to~4. An example of what these tell us is the following corollary about fully open CRNs:
\begin{prop}
\label{propMAfo}
If a fully open CRN $\mathcal{R}$ with mass action kinetics admits nondegenerate (resp., stable) oscillation, then so does any fully open CRN with mass action kinetics which includes $\mathcal{R}$ as an induced subnetwork.
\end{prop}

The definitions required to make this result precise will follow. Proposition~\ref{propMAfo} is the specialisation for mass action kinetics of a result with more general kinetic assumptions, Proposition~\ref{coropeninduced}, (see Remark~\ref{remMAfo}) which is a natural starting point for some computational exploration on small fully open CRNs admitting oscillation. It is worth noting at the outset that Proposition~\ref{propMAfo} fails if the CRNs are not assumed to be fully open. An example is provided in the concluding section (Example~\ref{exinherit}). 

Much of the mathematical literature on oscillation in CRNs has focussed on conditions which forbid oscillation, or forbid stable oscillation of the kind which might be observed in numerical simulations, or forbid bifurcations leading to oscillation. For CRNs with mass action kinetics, there are the original results of deficiency theory \cite{horn72,hornjackson,feinberg0,feinberg}; for CRNs with more general kinetics there are results based on the theory of monotone dynamical systems (\cite{banajidynsys,angelileenheersontag,donnellbanaji,banajimierczynski} for example), and algebraic approaches (\cite{abphopf} for example). Various papers which do not directly treat CRNs also have natural applications to forbidding oscillation or stable oscillation in CRNs, including the work of Angeli, Hirsch and Sontag on ``coherent'' systems \cite{angelihirschsontag}, and of Li and Muldowney on generalised Bendixson's criteria \cite{li_muldowney_1993, li_muldowney_1996, li_muldowney_2000}. On the other hand oscillation has been shown to occur in numerical studies of various CRNs of interest (for example, \cite{dicera, WolfOsci, Kholodenko.2000aa, Qiao.2007aa}). Aside from numerical work, there exists an important strand of theory drawing on approaches in convex and toric geometry which provides {\em sufficient} conditions for Hopf bifurcations in CRNs with mass action and generalised mass action kinetics \cite{eiswirth91, eiswirth96,gatermann, errami2015}. These approaches lead to algorithms for the determination of parameter regions where Hopf bifurcation occurs. Other papers treating the question of sufficient conditions for oscillation in chemical reaction networks include \cite{minchevaroussel} and \cite{domijan}. 

The work here is aimed at closing the gap between theory which forbids oscillation and examples of oscillatory networks or particular sufficient conditions for oscillation. It is likely that many examples of CRNs admitting oscillation in fact oscillate because they inherit this property from a smaller CRN which admits oscillation, and the goal is then to identify an appropriate notion of inheritance, and minimal oscillatory CRNs in some sense. The importance of inheritance approaches is increasingly recognised. In \cite{ConradiShiuPTM}, Conradi and Shiu pose a question closely related to the main question in this paper, namely whether Hopf bifurcation is preserved when CRNs are modified in natural ways. The problem of identifying a ``minimal'' oscillatory subnetwork was tackled for the biologically important MAPK cascade in \cite{hadac}.

Computational work on fully open CRNs towards the end of the paper confirms the practical usefulness of inheritance approaches. As oscillation may occur in very small regions of parameter space, it may be hard to find by brute force in numerical simulations, even where it is straightforward to predict its occurrence by inheritance results. Finding a single small oscillatory CRN on the other hand immediately gives us knowledge of a large number of CRNs which inherit this oscillation. Ultimately, the hope is that examining CRNs which can neither be proven to forbid oscillation nor be shown to oscillate (using numerics, known sufficient conditions for oscillation, or inheritance results such as here) may lead to new theorems about necessary conditions for oscillation.

\subsection{Notational preliminaries}

\begin{notation}[Nonnegative and positive vectors]
\label{notpos}
A real vector $x = (x_1, \ldots, x_n)^{\mathrm{t}}$ is nonnegative (resp., positive) if $x_i\geq 0$ (resp., $x_i > 0$) for each $i$, and we refer to the nonnegative (resp., positive) orthant in $\mathbb{R}^n$ as $\mathbb{R}^n_{\geq 0}$ (resp., $\mathbb{R}^n_{\gg 0}$). Subsets of $\mathbb{R}^n_{\gg 0}$ are referred to as positive.
\end{notation}

\begin{notation}[Vector of ones]
$\mathbf{1}$ denotes a vector of ones whose length is inferred from the context. 
\end{notation} 

\begin{notation}[Identity matrix]
$I_n$ is the $n \times n$ identity matrix.
\end{notation} 

\begin{notation}[Set theoretic inverse]
Given sets $X,Y$ and a function $f\colon X \to Y$, not necessarily invertible, $f^{-1}$ will generally refer to the set theoretic inverse, namely, given $Y_0 \subseteq Y$, $f^{-1}(Y_0) = \{x \in X\colon f(x) \in Y_0\}$.
\end{notation}

\begin{notation}[Monomials, vector of monomials]
\label{notmon}
Given $x=(x_1,\ldots, x_n)^{\mathrm{t}}$ and $a = (a_1,\ldots, a_n)$, $x^a$ is an abbreviation for the (generalised) monomial $\prod_ix_i^{a_i}$. If $A$ is an $m \times n$ matrix with rows $A_1, \ldots, A_m$, then $x^A$ means the vector of (generalised) monomials $(x^{A_1}, x^{A_2}, \ldots, x^{A_m})^{\mathrm{t}}$. 
\end{notation} 

\begin{notation}[Entrywise product]
\label{nothad}
Given two matrices $A$ and $B$ with the same dimensions, $A \circ B$ will refer to the entrywise (or Hadamard) product of $A$ and $B$, namely $(A\circ B)_{ij} = A_{ij}B_{ij}$. 
\end{notation}

\section{Periodic orbits} We remind the reader of some standard results from Floquet theory (Chapters 3 and 4 of \cite{HaleOsci} for example) as needed here. Let $X \subseteq \mathbb{R}^r$ be open, $F\colon X \to \mathbb{R}^r$ be $C^1$, and consider the ODE 
\begin{equation}
\label{Floq0}
\dot x = F(x)
\end{equation}
on $X$. Assume that (\ref{Floq0}) has a nontrivial periodic solution $\theta\colon \mathbb{R} \to X$ with smallest positive period $T$, and with corresponding periodic orbit $\mathcal{O}:=\mathrm{im}\,\theta$. The variational equation about $\theta$ is
\begin{equation}
\label{Floq1}
\dot z= DF(\theta(t))z.
\end{equation}
$DF(\theta(t))$ is an $r \times r$ $T$-periodic matrix and Floquet theory tells us that any fundamental matrix solution $Z(t)$ of (\ref{Floq1}) can be written in the form
\[
Z(t) = A(t)e^{tB}
\]
where $A$ is a nonsingular $T$-periodic matrix, and $B$ is a constant matrix. The eigenvalues of $e^{TB}$ are termed the {\em characteristic multipliers} (or {\em Floquet multipliers}) of $\mathcal{O}$. If $Z(0) = I$, then $A(T) = A(0) = I$, in which case the characteristic multipliers are the eigenvalues of $Z(T)$. $\mathcal{O}$ is termed {\em hyperbolic} (resp., {\em linearly stable}) if $r-1$ of its characteristic multipliers are disjoint from (resp., inside) the unit circle in $\mathbb{C}$.  Hyperbolicity (resp., linear stability) of a periodic orbit is precisely hyperbolicity (resp., linear stability) of the associated fixed point of any Poincar\'e map constructed on a section transverse to the periodic orbit: see Chapter~10 onwards of \cite{wiggins}, for example. Hyperbolic periodic orbits survive under sufficiently small perturbations of vector fields in a sense made precise in Lemma~\ref{lemreg} below. Linear stability of a periodic orbit implies asymptotic orbital stability, namely that forward trajectories of all sufficiently nearby initial conditions converge to the periodic orbit (Theorem~4.2 in \cite{HaleOsci}).

The following is a well-known result of regular perturbation theory. $d_{\mathrm{H}}(\cdot, \cdot)$ denotes the Hausdorff distance between nonempty compact subsets of Euclidean space. 
\begin{lemma1}
\label{lemreg}
Let $X \subseteq \mathbb{R}^r$ be open, $\epsilon' > 0$ and $F\colon X \times (-\epsilon', \epsilon')\to \mathbb{R}^r$ be $C^1$. Consider the $\epsilon$-dependent family of ODEs on $X$
\begin{equation}
\specialnumber{${}_\epsilon$}\label{eqnFloqe}
\dot x = F(x,\epsilon)\,.
\end{equation} 
Suppose that \specialeqref{eqnFloqe}{${}_0$} has a nontrivial hyperbolic (resp., linearly stable) $T$-periodic orbit $\mathcal{O} \subseteq X$. Then there exists $\epsilon_0>0$ s.t. for $\epsilon \in (-\epsilon_0, \epsilon_0)$ \specialeqref{eqnFloqe}{${}_\epsilon$} has a hyperbolic (resp., linearly stable) periodic orbit $\mathcal{O}_\epsilon$ satisfying $\lim_{\epsilon \to 0} d_{\mathrm{H}}(\mathcal{O}_\epsilon, \mathcal{O}) = 0$ and with period $T_\epsilon$ satisfying $\lim_{\epsilon \to 0}T_\epsilon = T$. 
\end{lemma1}
\begin{pf}
These claims are proved, for example, by constructing a family of Poincar\'e maps $\Pi_\epsilon$ for \specialeqref{eqnFloqe}{${}_\epsilon$} and applying the implicit function theorem at the fixed point of $\Pi_0$ corresponding to $\mathcal{O}$ as described in Section IV of \cite{Fenichel79}. \hfill$\square$ \end{pf}

We next consider some specialisations of Floquet theory to systems with linear first integrals relevant to the study of CRNs. Let $x \in \mathbb{R}^n_{\gg 0}$, $v\colon\mathbb{R}^n_{\gg 0} \to \mathbb{R}^m$ be $C^1$, $\Gamma$ be an $n \times m$ real matrix of rank $r$, and consider the ODE
\begin{equation}
\label{Floq2}
\dot x = \Gamma v(x).
\end{equation}
Assume that (\ref{Floq2}) has a nontrivial positive periodic orbit $\mathcal{O}$ (see Notation~\ref{notpos}), namely there exists some periodic solution $\theta\colon \mathbb{R} \to \mathbb{R}^n_{\gg 0}$ of (\ref{Floq2}) with smallest period $T>0$ and with $\mathcal{O}:=\mathrm{im}\,\theta$. Clearly, $S_{\mathcal{O}} := (\mathcal{O} + \mathrm{im}\,\Gamma) \cap \mathbb{R}^n_{\gg 0}$ is locally invariant under (\ref{Floq2}). If $r \neq n$, then $\mathcal{O}$ cannot be hyperbolic or linearly stable in the senses defined above. However, our interest is in whether it is hyperbolic (resp., linearly stable) {\em relative} to $S_{\mathcal{O}}$. Associated with $S_{\mathcal{O}}$ are $r$ characteristic multipliers and we would like to know whether $r-1$ of these are disjoint from (resp., inside) the unit circle. The single remaining multiplier associated with $S_{\mathcal{O}}$, corresponding to travel along the periodic orbit, is $1$, while the additional $n-r$ multipliers associated with directions transverse to $S_{\mathcal{O}}$ are also easily shown all to be $1$.

An explicit calculation of the multipliers of $\mathcal{O}$ relative to $S_{\mathcal{O}}$ is needed in certain proofs to follow. This proceeds as follows. Choose $x_0 \in S_{\mathcal{O}}$ and choose $\Gamma_0$ to be any matrix whose columns form a basis for $\mathrm{im}\,\Gamma$. Define $Q$ by $\Gamma = \Gamma_0Q$, and define the bijection $h\colon \mathbb{R}^r \to x_0 + \mathrm{im}\,\Gamma$ by $h(y) = x_0 + \Gamma_0 y$. Note that $W:= h^{-1}(S_{\mathcal{O}})$ is an open subset of $\mathbb{R}^r$, and $\left.h\right|_{W}$ is an affine bijection between $W$ and $S_{\mathcal{O}}$. Setting $x=h(y)$ we get, for the evolution of $y$:
\begin{equation}
\label{Floq3}
\dot y = Qv(x_0 + \Gamma_0y)\,.
\end{equation}
(\ref{Floq3}) has a $T$-periodic solution $\psi\colon \mathbb{R} \to W$ defined by $\psi(t) = h^{-1}(\theta(t))$. Let $\mathcal{O}':=\mathrm{im}\,\psi = h^{-1}(\mathcal{O})$ be the corresponding periodic orbit. The multipliers of $\mathcal{O}'$ are precisely the multipliers of $\mathcal{O}$ relative to $S_{\mathcal{O}}$. By definition $\mathcal{O}'$ is hyperbolic (resp., linearly stable) if it has $r-1$ characteristic multipliers disjoint from (resp., inside) the unit circle in $\mathbb{C}$. This motivates the following definitions: 

\begin{def1}[NPPO, SPPO]
\label{defNPPO}
Let $\mathcal{O}$ be a positive periodic orbit of (\ref{Floq2}). With $h$ defined as above, $\mathcal{O}$ is a {\em nondegenerate positive periodic orbit (NPPO)} of (\ref{Floq2}) if $\mathcal{O}':= h^{-1}(\mathcal{O})$ is a hyperbolic periodic orbit of (\ref{Floq3}). $\mathcal{O}$ is a {\em linearly stable positive periodic orbit (SPPO)} of (\ref{Floq2}) if $\mathcal{O}':= h^{-1}(\mathcal{O})$ is a linearly stable periodic orbit of (\ref{Floq3}). An SPPO is clearly also an NPPO.
\end{def1}

\begin{figure}[h]
\begin{center}
\begin{tikzpicture}[scale=0.9]


\node at (2.3,2.3) {$\mathbb{R}^r$};

\draw [->, line width=0.04cm] (-3,0) -- (2.5,0);
\draw [->, line width=0.04cm] (0,-2.5) -- (0,2.5);
\fill[color=black!30, fill opacity=0.7] (-1,2) -- (2.5,-1) -- (-2.5,-2) -- cycle;

\begin{scope}[scale=0.7,cm={1.2,0.1,0.3,1,(-7.3cm,-4.3cm)}]
\draw[->, line width=0.04cm] (5.5, 2) .. controls (6, 2) and (6,3) .. (6,3.5);
\draw[->, line width=0.04cm] (6, 3.5) .. controls (6, 4) and (4,5) .. (3.5,4.5);
\draw[->, line width=0.04cm] (3.5, 4.5) .. controls (3, 4) and (5,2) .. (5.5,2);
\end{scope}
\node at (0.25,1.4) {$W=h^{-1}(S_\mathcal{O})$};
\node at (2.25,-0.9) {$\mathcal{O}' = \mathrm{im}\,\psi = h^{-1}(\mathcal{O})$};

\node at (4,1.6) {$h$};
\draw[->, line width=0.04cm] (2.5, 1) .. controls (3.5, 1.3) and (4.5,1.3) .. (5.5,1);

\begin{scope}[scale=0.6,xshift=10cm, yshift=-4cm]

\node at (9,8) {$\mathbb{R}^n$};

\draw[->, line width=0.04cm] (-0.5,0) -- (10,0);
\draw[->, line width=0.04cm] (0,-0.5) -- (0,8);
\draw[->, line width=0.04cm] (-0.5,-0.25) -- (8,4);

\fill[color=black!30, fill opacity=0.7] (0,7) -- (7,3.5) -- (6.5,0) -- cycle;

\draw[->, line width=0.04cm] (5.5, 2) .. controls (6, 2) and (6,3) .. (6,3.5);
\draw[->, line width=0.04cm] (6, 3.5) .. controls (6, 4) and (4,5) .. (3.5,4.5);
\draw[->, line width=0.04cm] (3.5, 4.5) .. controls (3, 4) and (5,2) .. (5.5,2);
\fill (5.1,3.7) circle (3pt);
\node at (4.85,3.35) {$x_0$};
\node at (7.3,2.1) {$\mathcal{O} = \mathrm{im}\,\theta$};
\node at (2.1,5.5) {$S_\mathcal{O}$};

\end{scope}

\end{tikzpicture}
\end{center}
\caption{\label{figschematic} $h$ defines an affine embedding of $\mathbb{R}^r$ into $\mathbb{R}^n$, illustrated in the case $r=2$ and $n=3$. The image of $h$ is $x_0 + \mathrm{im}\,\Gamma$, assumed to include a positive periodic orbit $\mathcal{O}$, and $h$ thus defines local coordinates on $S_\mathcal{O}$, the positive stoichiometry class of $\mathcal{O}$. Of interest is the hyperbolicity or linear stability of $\mathcal{O}$ relative to $S_\mathcal{O}$, and by definition $\mathcal{O}$ is an NPPO (resp., SPPO) if $\mathcal{O}' = h^{-1}(\mathcal{O})$ is nondegenerate (resp., linear stable).}
\end{figure}
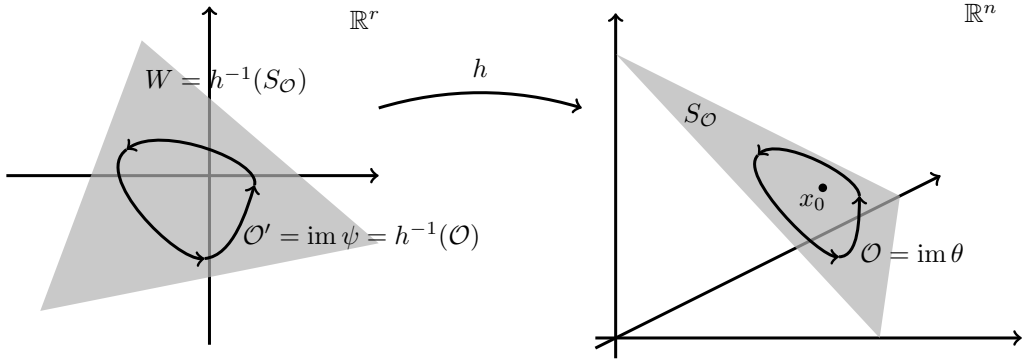

The use of the transformation $h$ to define new coordinates on $S_\mathcal{O}$ is illustrated schematically in Figure~\ref{figschematic}.

\begin{remark}
Note that the overloading of the term ``linearly stable'' in Definition~\ref{defNPPO} is an abuse of terminology which should cause no confusion: if $\mathrm{im}\,\Gamma = \mathbb{R}^n$, then linearly stable has its usual meaning; if $\mathrm{im}\,\Gamma \neq \mathbb{R}^n$, then no periodic orbit of (\ref{Floq2}) can truly be linearly stable, and linear stability is taken to mean linear stability relative $\mathrm{im}\,\Gamma$.
\end{remark}

We can easily verify that Definition~\ref{defNPPO} makes sense: different choices of $x_0$ or $\Gamma_0$ lead to the same characteristic multipliers. To see this, recall that according to Floquet theory the variational equation of (\ref{Floq3}) about $\psi(t) = h^{-1}(\theta(t))$, namely,
\begin{equation}
\label{Floq3a}
\dot z= QDv(\theta(t))\Gamma_0 z
\end{equation}
has a fundamental matrix solution $Z(t)$ which can be written $Z(t) = A(t)e^{tB}$ with $A$ a nonsingular $T$-periodic matrix and $B$ a constant matrix. The characteristic multipliers associated with $\psi$ are the eigenvalues of $e^{TB}$. Now suppose we make some different choices $x_0' \in S_{\mathcal{O}}$ and $\Gamma_0'$ and let $h'\colon \mathbb{R}^r \to S_{\mathcal{O}}$ be defined by $h'(y) := x_0'+\Gamma_0' y$. As the columns of $\Gamma_0'$ are a basis for $\mathrm{im}\,\Gamma$, $\Gamma_0' = \Gamma_0R$ where $R$ is a nonsingular $r \times r$ matrix. Thus $\Gamma = \Gamma_0Q = \Gamma_0'R^{-1}Q$. With $x = h'(y)$, we get the evolution on $W':=h'^{-1}(S_{\mathcal{O}})$
\begin{equation}
\label{Floq4}
\dot y = R^{-1}Qv(x_0' + \Gamma_0'y)
\end{equation}
with $T$-periodic solution $\psi'\colon \mathbb{R} \to W'$ defined by $\psi'(t) = h'^{-1}(\theta(t))$. The variational equation of (\ref{Floq4}) about $\psi'(t)$ is
\begin{equation}
\label{Floq4a}
\dot z= R^{-1}[QDv(\theta(t))\Gamma_0]R z. 
\end{equation}
Then $Z'(t):=R^{-1}Z(t) = R^{-1}A(t)e^{tB}$ is a fundamental matrix solution of (\ref{Floq4a}) with $R^{-1}A(t)$ clearly a $T$-periodic matrix. Thus the characteristic multipliers associated with the solution $\psi'(t)$ of (\ref{Floq4}) are again the eigenvalues of $e^{TB}$, i.e., those associated with the solution $\psi(t)$ of (\ref{Floq3}).

\section{Background on CRNs}
\label{secbackground}

As the framework and terminology closely follow that of \cite{banajipantea}, the reader is referred to this paper for some of the detail. The goal is to remain precise while minimising the extensive preamble on basic notation, terminology and definitions which accompanies many papers on CRNs. We consider a CRN involving $n$ chemical species $X_1, \ldots, X_n$.

\begin{def1}[Complexes, the zero complex, stoichiometry]A {\em complex} is a formal linear combination of species. If $a = (a_1, \ldots, a_n)^{\mathrm{t}}$ is a nonnegative integer vector, then $a\cdot X := a_1X_1 + a_2X_2 + \cdots + a_nX_n$ is a complex. $a_i$ is the {\em stoichiometry} of $X_i$ in the complex $a \cdot X$. The {\em zero complex} $0X_1 + \cdots + 0X_n$ is denoted $0$. 
\end{def1}

An irreversible reaction is an ordered pair of complexes, termed the {\em source complex} (or left hand side) and the {\em target complex} (or right hand side). We always assume that the source and target complexes are distinct. A reversible reaction may be considered either as two irreversible reactions or, equivalently, as an unordered pair of (distinct) complexes. A CRN is a set of species and a set of reactions. We adopt the common convention that the reactions of a CRN are distinct. However, for technical reasons, we do not forbid {\em a priori} the possibility that some chemical species occurs in a CRN but participates in none of its reactions. 

\begin{def1}[Flow reaction, fully open CRN, fully open extension of a CRN]
For the purposes of this paper, reactions of the form $0 \rightarrow A$ or $A \rightarrow 0$ are referred to as {\em flow reactions}, while all others are {\em non-flow reactions} (even where these clearly violate any conservation laws: for example $2A \rightarrow 0$ or $0 \rightarrow A+B$ are referred to as non-flow reactions). A CRN involving species $X_1, \ldots, X_n$ is {\em fully open} if it includes all flow reactions $0 \rightleftharpoons X_i$ ($i = 1, \ldots, n$). Note that if, for example, a CRN includes all reactions of the form $0 \rightleftharpoons 2X_i$, but not all the reactions $0 \rightleftharpoons X_i$, we do not refer to it as fully open. The {\em fully open extension} of a CRN $\mathcal{R}$ is the smallest fully open CRN containing all the reactions of $\mathcal{R}$, namely the CRN created by adjoining to $\mathcal{R}$ any flow reactions which are absent from $\mathcal{R}$. 
\end{def1}

\subsection{Combinatorial representations of CRNs}
CRNs are combinatorial objects which give rise to dynamical systems in different ways depending on various modelling choices. The most common combinatorial representation of a CRN is via its {\em complex graph} \cite{horn}, a digraph whose vertices are complexes and whose arcs correspond to (irreversible) reactions. For example, the reaction $X_1+2X_2 \rightarrow X_3$ is an ordered pair of complexes naturally represented as an arc from source complex $X_1+2X_2$ to target complex $X_3$. The set of species and the complex graph together make up a formal description of the CRN.

An alternative representation, particularly useful when discussing isomorphism of CRNs, is a {\em Petri net (PN) graph} \cite{angelipetrinet}, an edge-weighted bipartite digraph, defined in the form used here in \cite{banajipanteaMPNE}. The PN graph of a CRN $\mathcal{R}$, denoted $PN(\mathcal{R})$, has two vertex sets $V_S$ (species vertices) and $V_R$ (reaction vertices) identified with the species and the reactions of $\mathcal{R}$. Given $X_i \in V_S$ and $R_j \in V_R$, there exists an arc $X_iR_j$ (resp., $R_jX_i$) with weight $w$ if and only if the species corresponding to $X_i$ occurs with stoichiometry $w>0$ in the source complex (resp., target complex) of the reaction corresponding to $R_j$. Arc weights of $1$ are omitted from drawings for neatness. An unlabelled PN graph is referred to as a {\em motif}.

CRNs $\mathcal{R}_1$ and $\mathcal{R}_2$ are isomorphic if $PN(\mathcal{R}_1)$ and $PN(\mathcal{R}_2)$ are isomorphic in a natural sense, namely there exists a relabelling of the vertices of $PN(\mathcal{R}_1)$ which preserves the bipartition and gives $PN(\mathcal{R}_2)$. Given CRNs $\mathcal{R}_1$ and $\mathcal{R}_2$, we say that $\mathcal{R}_1$ is an {\em induced subnetwork} of $\mathcal{R}_2$, and write $\mathcal{R}_1 \leq \mathcal{R}_2$, if $PN(\mathcal{R}_1)$ is a vertex-induced subgraph of $PN(\mathcal{R}_2)$. Clearly, the induced subnetwork relationship induces a partial order on the set of CRNs as discussed in \cite{banajiCRNcount,banajipanteaMPNE}. Note that if $\mathcal{R}_1 \leq \mathcal{R}_2$, the occurrence of a reaction $R$ in both $\mathcal{R}_1$ and $\mathcal{R}_2$ does {\em not} mean that $R$ is, physically speaking, the same reaction with the same source and target complexes in $\mathcal{R}_2$ as in $\mathcal{R}_1$: identifying reactions with (labelled) vertices in a PN graph means that they maintain their identity as graph theoretic modifications are carried out equivalent to inserting or deleting species. If $\mathcal{R}_1 \leq \mathcal{R}_2$, and both have the same set of species, we say that $\mathcal{R}_1$ is a {\em reaction-induced subnetwork} of $\mathcal{R}_2$, and write $\mathcal{R}_1 \leq_{R} \mathcal{R}_2$. If $\mathcal{R}_1 \leq \mathcal{R}_2$, and both have the same set of reactions, we say that $\mathcal{R}_1$ is a {\em species-induced subnetwork} of $\mathcal{R}_2$, and write $\mathcal{R}_1 \leq_{S} \mathcal{R}_2$. Some of the definitions are illustrated in the following example. 

\begin{example}
\label{examplePN}
Consider the following CRN:
\[
X+Y \rightarrow 2Y,\quad Y+Z \rightarrow X \rightarrow W+Z, \quad W \rightarrow X. \tag{$\mathcal{R}$}
\]
$\mathcal{R}$ involves $4$ species $\{W,X,Y,Z\}$, $6$ complexes $\{W, X, 2Y,X+Y, Y+Z, W+Z\}$ and four (irreversible) reactions. The complex graph of $\mathcal{R}$ is shown below to the left and the PN graph in the centre. Removing the highlighted vertices and their incident arcs leads to the induced subnetwork
\[
X+Y \rightarrow 2Y, \quad Y+Z \rightarrow X, \tag{$\mathcal{R}_1$}
\]
represented in unlabelled form with species vertices as open circles and reaction vertices as filled circles to the right. Note that two reactions and a species were removed from $\mathcal{R}$ to obtain $\mathcal{R}_1$, and so the subnetwork is neither species-induced nor reaction-induced. 
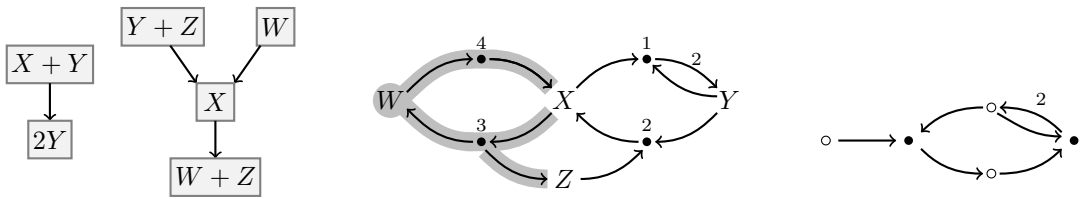
\begin{center}
\begin{tikzpicture}[scale=1,transition/.style={rectangle,draw=black!50,fill=black!5,thick,inner sep=0pt,minimum size=5mm}]
\node[transition] at (0,1.5) {$\,X+Y\,$};
\draw [->, thick] (0,1.25) --(0,0.75);
\node[transition] at (0,0.5) {$\,2Y\,$};

\node[transition] at (1.5,2) {$\,Y+Z\,$};
\node[transition] at (3,2) {$\,W\,$};
\draw [->, thick] (1.6,1.75) --(1.95,1.25);
\draw [->, thick] (2.8,1.75) --(2.45,1.25);
\node[transition] at (2.2,1) {$\,X\,$};
\draw [->, thick] (2.2,0.75) --(2.2,0.25);
\node[transition] at (2.2,0) {$\,W+Z\,$};

\end{tikzpicture}
\hspace{0.7cm}
\begin{tikzpicture}[scale=1.1]

%
%
%
%
\draw [-,color=black!25, line width=0.25cm] (0,0) .. controls (0.3,0.3) and (0.5,0.5) .. (0.9,0.5);
\draw [-,color=black!25, line width=0.25cm] (1.85,0.15) .. controls (1.7,0.3) and (1.5,0.5) .. (1.1,0.5);
\draw [-,color=black!25, line width=0.25cm] (0.8,0.5) -- (1.2,0.5);

\draw [-,color=black!25, line width=0.25cm] (0,0) .. controls (0.3,-0.3) and (0.5,-0.5) .. (0.9,-0.5);
\draw [-,color=black!25, line width=0.25cm] (1.85,-0.15) .. controls (1.7,-0.3) and (1.5,-0.5) .. (1.1,-0.5);
\draw [-,color=black!25, line width=0.25cm] (0.8,-0.5) -- (1.2,-0.5);

\draw [-,color=black!25, line width=0.25cm] (1.05,-0.6) .. controls (1.2,-0.75) and (1.4,-0.95) .. (1.8,-0.95);
\fill[color=black!25] (-0.1,0) circle (6pt);

\draw [->, thick] (1.05,-0.6) .. controls (1.2,-0.75) and (1.4,-0.95) .. (1.8,-0.95);
\draw [<-, thick] (2.95,-0.6) .. controls (2.8,-0.75) and (2.6,-0.95) .. (2.2,-0.95);
\node at (2,-0.95) {$Z$};
\draw [->, thick] (0.1,0.1) .. controls (0.3,0.3) and (0.5,0.5) .. (0.9,0.5);
\draw [<-, thick] (1.85,0.15) .. controls (1.7,0.3) and (1.5,0.5) .. (1.1,0.5);
\node at (1,0.7) {$\scriptstyle{4}$};
\fill[color=black] (1,0.5) circle (1.5pt);
\node at (1,-0.3) {$\scriptstyle{3}$};
\fill[color=black] (1,-0.5) circle (1.5pt);
\node at (3,0.7) {$\scriptstyle{1}$};
\fill[color=black] (3,0.5) circle (1.5pt);
\node at (3,-0.3) {$\scriptstyle{2}$};
\fill[color=black] (3,-0.5) circle (1.5pt);

\node at (-0.1,0) {$W$};
\node at (2,0) {$X$};
\node at (4,0) {$Y$};

\draw [<-, thick] (1.85,0.15) .. controls (1.7,0.3) and (1.5,0.5) .. (1.1,0.5);
\draw [->, thick] (1.85,-0.15) .. controls (1.7,-0.3) and (1.5,-0.5) .. (1.1,-0.5);
\draw [<-, thick] (0.1,-0.1) .. controls (0.3,-0.3) and (0.5,-0.5) .. (0.9,-0.5);
\draw [->, thick] (2.15,0.15) .. controls (2.3,0.3) and (2.5,0.5) .. (2.9,0.5);
\draw [<-, thick] (3.85,0.15) .. controls (3.7,0.3) and (3.5,0.5) .. (3.1,0.5);
\draw [<-, thick] (2.15,-0.15) .. controls (2.3,-0.3) and (2.5,-0.5) .. (2.9,-0.5);
\draw [->, thick] (3.85,-0.15) .. controls (3.7,-0.3) and (3.5,-0.5) .. (3.1,-0.5);
\draw[->, thick] (3.85,0.05) .. controls (3.5, 0.05) and (3.3, 0.2) .. (3.07,0.43);

\node at (3.6,0.5) {$\scriptstyle{2}$};
\end{tikzpicture}
\hspace{0.7cm}
\begin{tikzpicture}[scale=1.1]

\draw[color=black] (3,0.4) circle (1.5pt);
\draw[color=black] (3,-0.4) circle (1.5pt);
\draw[color=black] (1,0) circle (1.5pt);
\fill[color=black] (2,0) circle (1.5pt);
\fill[color=black] (4,0) circle (1.5pt);

\draw [->, thick] (1.15,0) --(1.85,0);
\draw [<-, thick] (2.15,0.1) .. controls (2.3,0.25) and (2.5,0.4) .. (2.9,0.4);
\draw [->, thick] (3.85,0.1) .. controls (3.7,0.25) and (3.5,0.4) .. (3.1,0.4);
\draw[<-, thick] (3.85,0.05) .. controls (3.5, 0.05) and (3.3, 0.2) .. (3.07,0.33);
\draw [->, thick] (2.15,-0.1) .. controls (2.3,-0.25) and (2.5,-0.4) .. (2.9,-0.4);
\draw [<-, thick] (3.85,-0.1) .. controls (3.7,-0.25) and (3.5,-0.4) .. (3.1,-0.4);

\node at (3.6,0.5) {$\scriptstyle{2}$};
\node at (2.4,-0.5) {$\textcolor{white}{\scriptstyle{2}}$};
\end{tikzpicture}
\end{center}
To preview the nature of results to follow, the motif on the right leads to stable periodic behaviour in fully open CRNs with mass action kinetics, and so the fully open extension of $\mathcal{R}$ with mass action kinetics admits an SPPO as a consequence of the presence of this motif. 

\end{example}

\subsection{ODE models of CRNs: basic definitions}We take the concentrations of chemical species to be nonnegative real numbers. Consider a CRN $\mathcal{R}$ involving $n$ chemical species $X_1, \ldots, X_n$ with corresponding concentration vector $x= (x_1, \ldots, x_n)^{\mathrm{t}}$, and $m$ irreversible reactions between the species. Orderings on the species and reactions are arbitrary but assumed fixed. Define nonnegative $n \times m$ matrices $\Gamma_l$ and $\Gamma_r$ as follows: $(\Gamma_l)_{ij}$ (resp., $(\Gamma_r)_{ij}$) is the stoichiometry of species $X_i$ on the left (resp., right) of reaction $j$. The {\em stoichiometric matrix} of $\mathcal{R}$ is $\Gamma=\Gamma_r-\Gamma_l$. The $j$th column of $\Gamma$ is termed the {\em reaction vector} for the $j$th reaction.

If the reactions of $\mathcal{R}$ proceed with rates $v_1(x), v_2(x),\ldots, v_m(x)$, we define the {\em rate function} of $\mathcal{R}$ to be $v(x) =(v_1(x), v_2(x),\ldots, v_m(x))^{\mathrm{t}}$. The evolution of the species concentrations is then governed by the ODE:
\begin{equation}
\label{genCRN}
\dot x = \Gamma v(x).
\end{equation}
If $v$ is defined and $C^1$ on $\mathbb{R}^n_{\gg 0}$ then (\ref{genCRN}) defines a local flow on $\mathbb{R}^n_{\gg 0}$, while if $v$ is defined and $C^1$ on $\mathbb{R}^n_{\geq 0}$ (namely, on an open subset of $\mathbb{R}^n$ containing $\mathbb{R}^n_{\geq 0}$) then, under physically reasonable assumptions on $v$ which ensure $\mathbb{R}^n_{\geq 0}$ is forward invariant, (\ref{genCRN}) defines a local semiflow on $\mathbb{R}^n_{\geq 0}$. See the introductory chapter of \cite{bhatiahajek} for definitions of local flows (there termed ``local dynamical systems'') and local semiflows (there termed ``local semi-dynamical systems''). 

$\mathrm{im}\,\Gamma$ is referred to as the {\em stoichiometric subspace} of the CRN. The nonempty intersection of a coset of $\mathrm{im}\,\Gamma$ with $\mathbb{R}^n_{\geq 0}$ (resp., $\mathbb{R}^n_{\gg 0}$) is a {\em stoichiometry class} (resp., {\em positive stoichiometry class}) of the CRN. If $\mathbb{R}^n_{\gg 0}$ (resp., $\mathbb{R}^n_{\geq 0}$) is forward invariant under the evolution defined by (\ref{genCRN}), then positive stoichiometry classes (resp., stoichiometry classes) are invariant under (\ref{genCRN}).

\subsection{Kinetics}
\label{seckin}
In order to state the results to follow with maximum applicability, we need some discussion of the rate functions of CRNs, namely the allowed functions $v$ in (\ref{genCRN}). The reader familiar with and primarily interested in mass action kinetics can skip directly to Proposition~\ref{propMAgen} below.

Given a CRN $\mathcal{R}$ with evolution governed by (\ref{genCRN}) we may assume that $v(x)$ belongs to some set of functions $\mathcal{K}$ with domain $\mathbb{R}^n_{\gg 0}$ and codomain $\mathbb{R}^m$. We refer to $\mathcal{K}$ as the {\em kinetics} of $\mathcal{R}$ and to the pair $(\mathcal{R},\mathcal{K})$ as a ``CRN with kinetics''. $\mathcal{K}$ may be finitely parameterised or a larger class of functions. Given a CRN with kinetics $(\mathcal{R}, \mathcal{K})$, and a given reaction $R$ in $\mathcal{R}$, the set of reaction rates for $R$ allowed by $\mathcal{K}$ is denoted $\mathcal{K}^{(R)}$.

When discussing kinetics it is assumed that a CRN consists of irreversible reactions: the allowed rates of a reversible reaction are derived by considering it as a pair of irreversible reactions. In each case below we assume that the CRN $\mathcal{R}$ involves $n$ species and $m$ (irreversible) reactions, and the $n \times m$ matrices $\Gamma_l$, $\Gamma_r$ and $\Gamma$ are defined as above. The following is a very large class of kinetics.

\begin{def1}[Positive general kinetics]
A rate function $v$ for $\mathcal{R}$ belongs to the class of {\em positive general kinetics} if and only if $v(x)$ is defined, positive-valued, and $C^1$ on $\mathbb{R}^n_{\gg 0}$ and satisfies for each $x \in \mathbb{R}^n_{\gg 0}$: (i) $\frac{\partial v_j}{\partial x_i} > 0$ if species $X_i$ occurs on the left of reaction $j$, (ii) $\frac{\partial v_j}{\partial x_i} = 0$ if species $X_i$ does not occur on the left of reaction $j$. Conditions (i) and (ii) can together be rephrased as ``the matrix $Dv(x)$ of partial derivatives of $v$ has the same sign pattern as $\Gamma_l$''. (Note that in \cite{banajipantea} condition (ii) was not spelled out explicitly, although it is implicit throughout.) 
\end{def1}

\begin{def1}[General kinetics]
A rate function $v$ for $\mathcal{R}$ belongs to the class of {\em general kinetics} if and only if $v(x)$ is defined and $C^1$ on $\mathbb{R}^n_{\geq 0}$, satisfies all the restrictions of positive general kinetics on $\mathbb{R}^n_{\gg 0}$, and $v_j(x) = 0$ if and only if $x_i = 0$ for some species $X_i$ occurring on the left of reaction $j$. $\mathbb{R}^n_{\geq 0}$ can easily be shown to be positively invariant for (\ref{genCRN}) under the assumption of general kinetics.
\end{def1}

\begin{def1}[Power-law kinetics, physical power-law kinetics, mass action kinetics]
A rate function $v$ for $\mathcal{R}$ belongs to the class of {\em power-law kinetics} if there exist $K \in \mathbb{R}^m_{\gg 0}$ and $M \in \mathbb{R}^{m \times n}$ such that $v(x) = K\circ x^M$ (recall Notation~\ref{notmon}~and~\ref{nothad} above). $K$ is the {\em vector of rate constants} and $M$ is the {\em matrix of exponents}. $v$ belongs to the class of {\em physical} power-law kinetics if, additionally, $M$ has the same sign pattern as $\Gamma_l^{\mathrm{t}}$, and of {\em mass action kinetics} if $M = \Gamma_l^{\mathrm{t}}$. 
\end{def1}

\begin{remark}[Fixed power-law kinetics]
Stating only that $\mathcal{R}$ has power-law kinetics, or physical power-law kinetics, implies that the entries of both $M$ and of $K$ are parameters which may vary. Stating that $\mathcal{R}$ has {\em fixed} power-law kinetics means that $M$ is fixed, while only the entries of $K$ are parameters which may vary. 
\end{remark}

\begin{remark}[Relationships between kinetic classes]
It is easily seen that physical power-law kinetics is a subclass of positive general kinetics, and that mass action kinetics is a particular case of fixed, physical power-law kinetics, and also of general kinetics. Further inclusions amongst classes of kinetics are detailed in \cite{banajipantea}.
\end{remark}

When we refer to $(\mathcal{R}, \mathcal{K})$ as a ``CRN with mass action kinetics'', or more briefly a ``mass action CRN'', this means that the set of allowed rate functions $\mathcal{K}$ is precisely that given by the assumption of mass action kinetics. A similar comment applies to other classes of kinetics. 

\begin{def1}[Derived power-law kinetics] 
\label{derived}
Let $(\mathcal{R}_1, \mathcal{K}_1)$ and $(\mathcal{R}_2, \mathcal{K}_2)$ be CRNs with fixed power-law kinetics and corresponding matrices of exponents $M_1$ and $M_2$. Let $\mathcal{R}_2$ have $n$ species and $m$ reactions and let $\mathcal{R}_1$ be an induced subnetwork of $\mathcal{R}_2$ with species indexed from $\alpha \subseteq \{1, \ldots, n\}$ and reactions indexed from $\beta \subseteq \{1, \ldots, m\}$. Then $\mathcal{K}_2$ is {\em derived} from $\mathcal{K}_1$ if $M_1 = M_2(\beta|\alpha)$, where $M_2(\beta|\alpha)$ is the submatrix of $M_2$ with rows from $\beta$ and columns from $\alpha$. 
\end{def1}

\begin{def1}[Scaling invariant kinetics] Let $(\mathcal{R}, \mathcal{K})$ be a CRN with kinetics. Then $\mathcal{K}$ is {\em scaling invariant} if, for each reaction $R$ of $\mathcal{R}$ and each $\epsilon > 0$, $F \in \mathcal{K}^{(R)}$ implies that $\epsilon F \in \mathcal{K}^{(R)}$. 
\end{def1}

\begin{remark}[Scaling invariant kinetics]
A CRN with any reasonable kinetics, including positive general kinetics, power-law kinetics, physical power-law kinetics, or any fixed power-law kinetics (including mass action) has kinetics which is scaling invariant: if $v_j(x)$ is an allowed reaction rate from one of these classes for reaction $j$, then so is $\epsilon v_j(x)$ for each $\epsilon > 0$. 
\end{remark}

\subsection{Extending the kinetics of an induced subnetwork} Consider CRNs $\mathcal{R}_1 \leq \mathcal{R}_2$ with $\mathcal{R}_1$ given kinetics $\mathcal{K}_1$. Are there natural ways of ``extending'' $\mathcal{K}_1$ to a kinetics $\mathcal{K}_2$ for $\mathcal{R}_2$? For example, it is reasonable and often mathematically convenient to assume that: 
\begin{itemize}[align=left,leftmargin=*]
\item Where a reaction of $\mathcal{R}_1$ occurs with the same source and target complexes in $\mathcal{R}_2$, the rates for this reaction allowed by $\mathcal{K}_2$ should include those allowed by $\mathcal{K}_1$. 
\item Where a reaction of $\mathcal{R}_1$ occurs in $\mathcal{R}_2$ with some new species involved, fixing the concentrations of the new species at some positive values should give back (at least) all the rate functions allowed by $\mathcal{K}_1$. 
\end{itemize}
These notions are formalised in the following two definitions. 

\begin{def1}[Reaction-extensions] 
\label{defreacext}
Consider CRNs with kinetics $(\mathcal{R}_1, \mathcal{K}_1)$ and $(\mathcal{R}_2, \mathcal{K}_2)$. Then $(\mathcal{R}_2, \mathcal{K}_2)$ is a {\em reaction-extension} of $(\mathcal{R}_1, \mathcal{K}_1)$, written $(\mathcal{R}_1, \mathcal{K}_1)\leq_{R} (\mathcal{R}_2, \mathcal{K}_2)$, if $\mathcal{R}_1$ is a reaction-induced subnetwork of $\mathcal{R}_2$ and $\mathcal{K}_1^{(R)} \subseteq\mathcal{K}_2^{(R)}$ for each reaction $R$ occurring in both $\mathcal{R}_1$ and $\mathcal{R}_2$. In other words, reactions which $\mathcal{R}_2$ inherits from $\mathcal{R}_1$ are allowed (at least) all the rate functions allowed by $\mathcal{K}_1$. 
\end{def1}

\begin{def1}[Species-extensions] 
\label{defspecext}
Consider CRNs with kinetics $(\mathcal{R}_1, \mathcal{K}_1)$ and $(\mathcal{R}_2, \mathcal{K}_2)$ and suppose that $\mathcal{R}_1$ is a species-induced subnetwork of $\mathcal{R}_2$. Let $\mathcal{R}_2$ have $n_2$ species $X_1, \ldots, X_{n_2}$ and assume, without loss of generality, that the species of $\mathcal{R}_1$ are $X_1, \ldots, X_{n_1}$ where $n_1 \leq n_2$. Let $\hat{x} = (x_1, \ldots, x_{n_1})$ and $\doublehat{x} = (x_{n_1 + 1}, \ldots, x_{n_2})$. Then $(\mathcal{R}_2, \mathcal{K}_2)$ is a {\em species-extension} of $(\mathcal{R}_1, \mathcal{K}_1)$, written $(\mathcal{R}_1, \mathcal{K}_1)\leq_{S} (\mathcal{R}_2, \mathcal{K}_2)$, if for each $v(\hat{x}) \in \mathcal{K}_1$ there exists $w(\hat{x},\doublehat{x}) \in \mathcal{K}_2$ such that $w(\hat{x},\mathbf{1}) = v(\hat{x})$, and such that if $v$ is $C^k$ on $\mathbb{R}^{n_1}_{\gg 0}$ (resp., $\mathbb{R}^{n_1}_{\geq 0}$), then $w$ is $C^k$ on $\mathbb{R}^{n_2}_{\gg 0}$ (resp., $\mathbb{R}^{n_2}_{\geq 0}$). (In the case $n_1=n_2$ and $\mathcal{K}_1 \subseteq \mathcal{K}_2$ we take $(\mathcal{R}_1, \mathcal{K}_1)\leq_{S} (\mathcal{R}_2, \mathcal{K}_2)$ to be trivially true.)
\end{def1}

\begin{lemma1}[Species-extensions]
\label{lemspecext}
Let $\mathcal{R}_1\leq_{S} \mathcal{R}_2$. Then the CRNs with kinetics $(\mathcal{R}_1, \mathcal{K}_1)$ and $(\mathcal{R}_2, \mathcal{K}_2)$ satisfy $(\mathcal{R}_1, \mathcal{K}_1) \leq_S (\mathcal{R}_2, \mathcal{K}_2)$ if any of the following hold:
\begin{enumerate}[align=left,leftmargin=*]
\item $\mathcal{K}_1$ and $\mathcal{K}_2$ are both given by positive general kinetics.
\item $\mathcal{K}_1$ and $\mathcal{K}_2$ are both given by power-law kinetics.
\item $\mathcal{K}_1$ and $\mathcal{K}_2$ are both given by physical power-law kinetics.
\item $\mathcal{K}_1$ and $\mathcal{K}_2$ are both given by mass action kinetics. 
\item $\mathcal{K}_1$ and $\mathcal{K}_2$ are both given by fixed power-law kinetics with $\mathcal{K}_2$ derived from $\mathcal{K}_1$ (see Definition~\ref{derived}).
\end{enumerate}
\end{lemma1}
\begin{pf}
Using the notation in Definition~\ref{defspecext}, for each rate function $v \in \mathcal{K}_1$, set $w(\hat{x}, \doublehat{x}) = v(\hat{x})\circ\doublehat{x}^{M}$ where, in cases (1)~to~(4), $M$ consists of the final $n_2-n_1$ columns of $\Gamma_l^{\mathrm{t}}$,
while in case (5) $M$ consists of the final $n_2-n_1$ columns of $M_2$, the matrix of exponents of $\mathcal{R}_2$.
In each case it is easily seen that $w \in \mathcal{K}_2$, that $w(\hat{x},\mathbf{1}) = v(\hat{x})$, and that if $v$ is $C^k$ on $\mathbb{R}^{n_1}_{\gg 0}$ (resp., $\mathbb{R}^{n_1}_{\geq 0}$ in case 4), then $w$ is $C^k$ on $\mathbb{R}^{n_2}_{\gg 0}$ (resp., $\mathbb{R}^{n_2}_{\geq 0}$ in case 4).
\end{pf}

\begin{def1}[Species-reaction-extensions] 
\label{defspecreacext}
Let $\mathcal{R}_1 \leq \mathcal{R}_2$ and consider CRNs with kinetics $(\mathcal{R}_1, \mathcal{K}_1)$ and $(\mathcal{R}_2, \mathcal{K}_2)$. Observe that there is a uniquely defined CRN $\mathcal{R}'$ satisfying $\mathcal{R}_1 \leq_S \mathcal{R}' \leq_{R} \mathcal{R}_2$. ($\mathcal{R}'$ is obtained by inserting all missing species into the reactions of $\mathcal{R}_1$, but without adding any new reactions.) Then $(\mathcal{R}_2, \mathcal{K}_2)$ is a {\em species-reaction-extension} of $(\mathcal{R}_1, \mathcal{K}_1)$ if there exists $\mathcal{K}'$ such that 
\[
(\mathcal{R}_1, \mathcal{K}_1) \leq_S (\mathcal{R}', \mathcal{K}') \leq_R (\mathcal{R}_2, \mathcal{K}_2).
\]
Intuitively, we first add in missing species and extend the kinetics of any modified reactions consistent with the species-extension condition, and then add in any remaining missing reactions. 
\end{def1}

\section{Results on the inheritance of NPPOs and SPPOs}
\label{secthms}
A CRN with kinetics $(\mathcal{R}, \mathcal{K})$ {\em admits} an NPPO (resp., SPPO) if there exists {\em some} rate function $v \in \mathcal{K}$ s.t. the associated ODE system (\ref{genCRN}) has an NPPO (resp., SPPO). A broad question is when, given CRNs with kinetics $(\mathcal{R}_1, \mathcal{K}_1)$ and $(\mathcal{R}_2, \mathcal{K}_2)$ related in some natural way, knowledge that one admits an NPPO (resp., SPPO) allows us to predict the same for the other. Four ``inheritance'' theorems in this direction will be proved below under varying kinetic assumptions. For the reader primarily interested in mass action kinetics, these can be summarised in a single corollary:
\begin{prop}
\label{propMAgen}
Let $\mathcal{R}$ and $\mathcal{R}'$ be CRNs, and suppose that $\mathcal{R}$ admits an NPPO (resp., SPPO) with mass action kinetics. Suppose that we create $\mathcal{R}'$ from $\mathcal{R}$ by
\begin{enumerate}[align=left,leftmargin=*]
\item adding to $\mathcal{R}$ a new reaction with reaction vector in the span of reaction vectors of $\mathcal{R}$; or
\item taking the fully open extension of $\mathcal{R}$; or
\item adding into some reactions of $\mathcal{R}$ a new species $Y$ which occurs with the same stoichiometry on both sides of each reaction in which it participates; or
\item adding into reactions of $\mathcal{R}$ a new species $Y$ in any way, while also adding the new reaction $0 \rightleftharpoons Y$. 
\end{enumerate}
Then, with mass action kinetics, $\mathcal{R}'$ admits an NPPO (resp., SPPO). 
\end{prop}
\begin{pf}
Claims 1 to 4 are immediate corollaries of Theorems~\ref{thmnewdepreac}~to~\ref{thmnewwithopen} below and the surrounding remarks. In order to apply the results we need only note that mass action kinetics is polynomial and hence certainly $C^2$, is scaling invariant, and that the assumptions imply that $\mathcal{R}'$ with mass action kinetics is a species-reaction extension of $\mathcal{R}$ with mass action kinetics. \hfill$\square$
\end{pf} 

Theorems~\ref{thmnewdepreac}~and~\ref{thmopenextension} require only basic regular perturbation theory to prove: in Theorem~\ref{thmnewdepreac} the application is almost trivial while in Theorem~\ref{thmopenextension} it takes a little more work to set up the problem. Theorem~\ref{thmtrivial} requires essentially no machinery to prove: the proof is almost immediate from the definitions. Theorem~\ref{thmnewwithopen} requires some results from singular perturbation theory. Theorems~\ref{thmnewdepreac}~and~\ref{thmnewwithopen} together imply an important corollary about fully open networks spelled out as Proposition~\ref{coropeninduced}.

In each of the following theorems, $\mathcal{R}$ is a CRN with $m$ reactions involving $n$ species $X_1, \ldots, X_n$ with concentrations $x_1, \ldots, x_n$. $\Gamma$, the stoichiometric matrix of $\mathcal{R}$, has rank $r$, $\Gamma_0$ is a matrix whose columns are a basis for $S := \mathrm{im}\,\Gamma$, and $Q$ is defined by $\Gamma = \Gamma_0Q$. Given a periodic orbit $\mathcal{O}$, $S_\mathcal{O}:=(\mathcal{O} + S) \cap \mathbb{R}^n_{\gg 0}$ is the positive stoichiometry class of $\mathcal{O}$, and $x_0$ is some point on $S_\mathcal{O}$ (recall Figure~\ref{figschematic}).

\begin{thm}[Adding a dependent reaction]
\label{thmnewdepreac}
Let $(\mathcal{R}, \mathcal{K})$ be a CRN with $C^1$ kinetics admitting an NPPO (resp., SPPO). Let $(\mathcal{R}', \mathcal{K}')$ be a reaction-extension of $(\mathcal{R}, \mathcal{K})$ created by adding to $\mathcal{R}$ a new irreversible reaction with $C^1$, scaling invariant, kinetics, and with reaction vector in the span of reaction vectors of $\mathcal{R}$. Then $(\mathcal{R}', \mathcal{K}')$ admits an NPPO (resp., SPPO). 
\end{thm}
\begin{pf}
Fix the rate function $v \in \mathcal{K}$ such that $\mathcal{R}$ has an NPPO (resp., SPPO) $\mathcal{O}$. Let the new reaction of $\mathcal{R}'$ be $a\cdot X \rightarrow a' \cdot X$. Define $\alpha = a'-a$ and define $c$ by $\alpha = \Gamma_0 c$. Consistent with the kinetic assumptions, set the rate of the new reaction to be $\epsilon f(x)$ where $f\colon \mathbb{R}^n_{\gg 0} \to \mathbb{R}$ is $C^1$ and $\epsilon$ is a parameter to be controlled (for example, with mass action kinetics the rate would be $\epsilon x^a$). The evolution of $\mathcal{R}'$ is then governed by:
\begin{equation}
\specialnumber{${}_\epsilon$}\label{eqCRN1}
\dot x = \Gamma v(x) + \epsilon \alpha f(x)  = \Gamma_0(Qv(x) +\epsilon c f(x)).
\end{equation}
Define $h \colon \mathbb{R}^r \to x_0+S$ by $h(z) := x_0+\Gamma_0 z$. $h$ is an affine bijection between $h^{-1}(S_\mathcal{O})$ and $S_\mathcal{O}$ and defines local coordinates on $S_\mathcal{O}$ via $x = h(z)$. $z$ evolves according to 
\begin{equation}
\specialnumber{${}_\epsilon$}\label{eqnThm1}
\dot z = Qv(x_0+\Gamma_0 z) + \epsilon c f(x_0+\Gamma_0 z).
\end{equation}
By definition, \specialeqref{eqnThm1}{${}_0$} has the hyperbolic (resp., linearly stable) periodic orbit $\mathcal{O}':=h^{-1}(\mathcal{O})$. By Lemma~\ref{lemreg} there exists $\epsilon_0>0$ s.t. for $\epsilon \in (-\epsilon_0, \epsilon_0)$ \specialeqref{eqnThm1}{${}_\epsilon$} has a hyperbolic (resp., linearly stable) periodic orbit $\mathcal{O}'_\epsilon$. Thus, for $\epsilon \in (-\epsilon_0, \epsilon_0)$, \specialeqref{eqCRN1}{${}_\epsilon$} has the NPPO (resp., SPPO) $\mathcal{O}_\epsilon := h(\mathcal{O}'_\epsilon)$. 
\hfill$\square$
\end{pf}

\begin{remark}[Adding the reverse of a reaction] 
\label{newdeprev}
Clearly, by Theorem~\ref{thmnewdepreac}, given a CRN $\mathcal{R}$ with kinetics from any $C^1$ class admitting an NPPO (resp., SPPO), adding the reverse of any existing reaction to $\mathcal{R}$ with $C^1$, scaling invariant, kinetics preserves this property. Thus if a CRN with, say, mass action kinetics admits an NPPO (resp., SPPO), then so does the corresponding reversible CRN with mass action kinetics. 
\end{remark}

\begin{remark}[Preservation of bifurcations when dependent reactions are added] 
In \cite{ConradiShiuPTM} Conradi and Shiu posed the question of whether Hopf bifurcations in CRNs are preserved when some irreversible reactions are made reversible. In fact, any generic bifurcation (\cite{wiggins} or \cite{kuznetsov} for example) survives the addition of dependent reactions with sufficiently smooth, scaling-invariant, kinetics. Although Theorem~\ref{thmnewdepreac} is not about bifurcation {\em per se}, the key idea in its proof is the construction of local coordinates on a stoichiometry class $S$ so that the vector field of $\mathcal{R}'$ in these local coordinates is a perturbation of the original vector field of $\mathcal{R}$. Suppose that some $C^r$ $k$-parameter family of vector fields $\mathcal{F}_\lambda$ on $S$ associated with $\mathcal{R}$ admits a nondegenerate codimension-$k$ bifurcation at $(x_0, \lambda_0)$. Then, as we see from \specialeqref{eqnThm1}{${}_\epsilon$}, addition of a new dependent reaction with $C^r$, scaling-invariant, kinetics gives rise, for each fixed $\epsilon$, to a new $C^r$, $k$-parameter, family $\mathcal{F}^\epsilon_\lambda$ of vector fields for $\mathcal{R}'$, $C^r$ close to $\mathcal{F}_\lambda$; for $r$ sufficiently large, and $\epsilon$ sufficiently small, the family $\mathcal{F}^\epsilon_\lambda$ will admit the same nondegenerate bifurcation. Analogous remarks apply to the other network modifications detailed in the theorems to follow. As a practical note, confirming that a given CRN does indeed admit a generic Hopf bifurcation at some parameter values is not always entirely straightfoward, as it may involve approximation of a parameter-dependent center manifold in order to confirm the nondegeneracy conditions.
\end{remark}

\begin{thm}[Adding inflows and outflows of all species]
\label{thmopenextension}
Let $(\mathcal{R}, \mathcal{K})$ be a CRN with $C^1$ kinetics admitting an NPPO (resp., SPPO). Suppose that $\mathcal{R}$ includes no flow reactions (i.e., no reactions of the form $0 \rightarrow X_i$ or $X_i \rightarrow 0$). Let $(\mathcal{R}', \mathcal{K}')$ be a reaction-extension of $(\mathcal{R}, \mathcal{K})$ created by adding to $\mathcal{R}$ all the reactions $0 \rightleftharpoons X_i$ ($i= 1, \ldots, n$) with kinetics from a class including mass action kinetics. Then $(\mathcal{R}', \mathcal{K}')$ admits an NPPO (resp., SPPO).
\end{thm}

\begin{pf}
Fix the rate function $v \in \mathcal{K}$ such that $\mathcal{R}$ has an NPPO (resp., SPPO) $\mathcal{O}$. Treat the $i$th inflow-outflow reaction as a single reversible reaction with mass action kinetics and forward and backwards rate constants $\epsilon (x_0)_i$ and $\epsilon$ respectively. The evolution of $\mathcal{R}'$ is then governed by:
\[
\dot x = \Gamma v(x) + \epsilon I_n (x_0 - x).
\]
Let $\Gamma_0' = [\Gamma_0|\Gamma_1]$ where $\Gamma_1$ is any $n \times (n-r)$ matrix chosen so that $\Gamma_0'$ has rank $n$. Observe that $\Gamma = \Gamma_0'\left(\begin{array}{c}Q\\0\end{array}\right)$. Define new coordinates $z = (\hat{z}, \doublehat{z}) \in \mathbb{R}^{r} \times \mathbb{R}^{n-r}$ by $x = h(z) := x_0 + \Gamma_0'z$. $h$ is an affine bijection between $W:=h^{-1}(\mathbb{R}^n_{\gg 0}) \subseteq \mathbb{R}^{r} \times \mathbb{R}^{n-r}$ and $\mathbb{R}^n_{\gg 0}$, and $z$ evolves according to 
\begin{equation}
\specialnumber{${}_\epsilon$}\label{eqnFO}
\frac{\mathrm{d}}{\mathrm{d}t}\left(\begin{array}{cc}\hat{z}\\ \doublehat{z}\end{array}\right) = \left(\begin{array}{c}Q\\0\end{array}\right)v(x_0+\Gamma_0' z) - \epsilon \left(\begin{array}{cc}\hat{z}\\ \doublehat{z}\end{array}\right)\,.
\end{equation}
Define $W_1:=W \cap (\mathbb{R}^r \times \{0\})$, so that $h(W_1) = S_{\mathcal{O}}$. Define $\mathcal{O}':=h^{-1}(\mathcal{O}) \subseteq W_1$ and define $\overline{\mathcal{O}} \subseteq \mathbb{R}^r$ by $\overline{\mathcal{O}} \times \{0\} = \mathcal{O}'$. $W_1$ is locally invariant for \specialeqref{eqnFO}{${}_\epsilon$}, and restricting \specialeqref{eqnFO}{${}_\epsilon$} to $W_1$ gives the differential equation:
\begin{equation}
\specialnumber{${}_\epsilon$}\label{eqnFOa}
\frac{\mathrm{d}\hat{z}}{\mathrm{d}t} = Qv(x_0+\Gamma_0 \hat{z}) - \epsilon \hat{z}\,. 
\end{equation}
By definition, $\mathcal{O}$ is an NPPO (resp., SPPO) of $\mathcal{R}$ if and only if $\overline{\mathcal{O}}$ is a hyperbolic (resp., linearly stable) periodic orbit of \specialeqref{eqnFOa}{${}_0$}. In this case, by Lemma~\ref{lemreg}, there exists $\epsilon_0 > 0$ s.t. for $\epsilon \in (-\epsilon_0, \epsilon_0)$, \specialeqref{eqnFOa}{${}_\epsilon$} has a hyperbolic (resp., linearly stable) periodic orbit $\overline{\mathcal{O}}_\epsilon$ close to $\overline{\mathcal{O}}$ with period $T_\epsilon$ close to $T$. It remains to show that $\mathcal{O}'_\epsilon:=\overline{\mathcal{O}}_\epsilon \times \{0\}$ is hyperbolic (resp., linearly stable) for \specialeqref{eqnFO}{${}_\epsilon$} for all sufficiently small $\epsilon > 0$. This will imply immediately that $\mathcal{O}_\epsilon := h(\mathcal{O}'_\epsilon)$ is an NPPO (resp., SPPO) of $\mathcal{R}'$. 

For each fixed $\epsilon \in (0, \epsilon_0)$, choose $\psi_\epsilon$ to be some solution of \specialeqref{eqnFOa}{${}_\epsilon$} with initial condition on $\overline{\mathcal{O}}_\epsilon$. The variational equation of \specialeqref{eqnFOa}{${}_\epsilon$} about $\psi_\epsilon$ is:
\begin{equation}
\specialnumber{${}_\epsilon$}\label{eqFOred}
\frac{\mathrm{d}\hat{\zeta}}{\mathrm{d}t}= [QDv(x_0+\Gamma_0 \psi_\epsilon(t))\Gamma_0 - \epsilon I_{r}]\hat{\zeta}\,.
\end{equation}
The fundamental matrix solution $\hat{Z}_\epsilon(t)$ of \specialeqref{eqFOred}{${}_\epsilon$} with $\hat{Z}_\epsilon(0) = I_r$ can be written $\hat{Z}_\epsilon(t) = A_\epsilon(t)e^{tB_\epsilon}$ where $A_\epsilon(t)$ is a nonsingular periodic matrix of period $T_\epsilon>0$ and $B_\epsilon$ is a constant matrix. Hyperbolicity (resp., linear stability) of $\overline{\mathcal{O}}_\epsilon$ for \specialeqref{eqnFOa}{${}_\epsilon$} means that $\hat{Z}_\epsilon(T_\epsilon) = e^{T_\epsilon B_\epsilon}$ has exactly one eigenvalue equal to $1$ with the remaining $r-1$ eigenvalues disjoint from (resp., inside) the unit circle.

For each $\psi_\epsilon$ chosen as above, $(\psi_\epsilon, 0)$ is clearly a periodic solution of \specialeqref{eqnFO}{${}_\epsilon$}, with image $\mathcal{O}'_\epsilon$. The full variational equation of \specialeqref{eqnFO}{${}_\epsilon$} about $(\psi_\epsilon, 0)$ is:
\begin{equation}
\label{eqnFOredb}
\frac{\mathrm{d}}{\mathrm{d}t}\left(\begin{array}{cc}\hat{\zeta}\\ \doublehat{\zeta}\end{array}\right) = \left(\begin{array}{cc}QDv(x_0+\Gamma_0 \psi_\epsilon(t))\Gamma_0 - \epsilon I_{r} & QDv(x_0+\Gamma_0 \psi_\epsilon(t))\Gamma_1\\0 & -\epsilon I_{n-r}\end{array}\right)\left(\begin{array}{cc}\hat{\zeta}\\ \doublehat{\zeta}\end{array}\right)\,.
\end{equation}
Our goal is to compute $Z_\epsilon(t)$, the fundamental matrix solution of (\ref{eqnFOredb}) satisfying $Z_\epsilon(0) = I$. Solving the second equation of (\ref{eqnFOredb})
gives $\doublehat{\zeta}(t) = e^{-\epsilon t}\doublehat{\zeta}(0)$. Substituting into the first equation of (\ref{eqnFOredb}) gives
\[
\frac{\mathrm{d}\hat{\zeta}}{\mathrm{d}t} = [QDv(x_0+\Gamma_0 \psi_\epsilon(t))\Gamma_0 - \epsilon I_{r}]\hat{\zeta} + e^{-\epsilon t}QDv(x_0+\Gamma_0 \psi_\epsilon(t))\Gamma_1\doublehat{\zeta}(0).
\]
Setting $\doublehat{\zeta}(0) = 0$ gives back \specialeqref{eqFOred}{${}_\epsilon$}. The above calculations give:
\[
Z_\epsilon(t) = \left(\begin{array}{cc}\hat{Z}_\epsilon(t) & A(t)\\0 & e^{-\epsilon t} I_{n-r}\end{array}\right),\quad \mbox{and hence,} \quad
Z_\epsilon(T_{\epsilon}) = \left(\begin{array}{cc}\hat{Z}_\epsilon(T_\epsilon) & A(T_\epsilon)\\0 & e^{-\epsilon T_\epsilon} I_{n-r}\end{array}\right).
\]
Here $A(t)$ is some matrix which can be determined by integration but which does not affect the subsequent argument. The characteristic multipliers of $\mathcal{O}'_\epsilon$ are precisely the eigenvalues of $\hat{Z}_\epsilon(T_\epsilon)$, namely the characteristic multipliers of $\overline{\mathcal{O}}_\epsilon$ for \specialeqref{eqnFOa}{${}_\epsilon$}, and the single value $e^{-\epsilon T_\epsilon}$ occurring with multiplicity $n-r$. As $\overline{\mathcal{O}}_\epsilon$ is a hyperbolic (resp., linearly stable) periodic orbit of \specialeqref{eqnFOa}{${}_\epsilon$}, and $e^{-\epsilon T_\epsilon}$ lies inside the unit circle for any $\epsilon > 0$, $T_\epsilon>0$, $\mathcal{O}'_\epsilon$ is a hyperbolic (resp., linearly stable) periodic orbit of \specialeqref{eqnFO}{${}_\epsilon$}, and consequently $\mathcal{O}_\epsilon= h(\mathcal{O}'_\epsilon)$ is an NPPO (resp., SPPO) of $\mathcal{R}'$.
\hfill$\square$
\end{pf}

\begin{remark}[Geometric interpretation of Theorem~\ref{thmopenextension}, the role of mass action kinetics]
Inflows and outflows were chosen to guarantee that $S_{\mathcal{O}}$ remained invariant for $\mathcal{R}'$: this necessitated mass action kinetics for the flow reactions. The construction ensured that for sufficiently small $\epsilon > 0$ $S_{\mathcal{O}}$ is exponentially attracting and the vector field of $\mathcal{R}'$ restricted to $S_{\mathcal{O}}$ is $\epsilon$-close to that of $\mathcal{R}$ restricted to $S_{\mathcal{O}}$, ensuring the existence on $S_{\mathcal{O}}$ of a hyperbolic (resp., linearly stable) periodic orbit $\mathcal{O}_\epsilon$ close to $\mathcal{O}$. 
\end{remark}

\begin{remark}[Theorem~\ref{thmopenextension} and fully open extensions]
Suppose that $(\mathcal{R},\mathcal{K})$ is any CRN with $C^1$ kinetics such that if $v$ is an allowed rate for some reaction $X_i \rightarrow 0$ of $\mathcal{R}$, then so is $v+\epsilon x_i$ for all sufficiently small $\epsilon > 0$, and if $v$ is an allowed rate for some reaction $0 \rightarrow X_i$ of $\mathcal{R}$, then so is $v+\epsilon$ for all sufficiently small $\epsilon > 0$. Then the condition that $\mathcal{R}$ excludes flow reactions can clearly be dropped in Theorem~\ref{thmopenextension}. In particular, if $(\mathcal{R},\mathcal{K})$ is any mass action CRN admitting an NPPO (resp., SPPO) then Theorem~\ref{thmopenextension} tells us that its fully open extension with mass action kinetics admits an NPPO (resp., SPPO). 
The same holds for CRNs with positive general kinetics. However, we cannot arrive at this conclusion for CRNs with arbitrary fixed physical power-law kinetics.
\end{remark}

\begin{thm}[Adding a trivial species]
\label{thmtrivial}
Let $(\mathcal{R}, \mathcal{K})$ be a CRN with $C^1$ kinetics admitting an NPPO (resp., SPPO). Let $(\mathcal{R}', \mathcal{K}')$ be a species-extension of $(\mathcal{R}, \mathcal{K})$ created by adding into some reactions of $\mathcal{R}$ a new species $Y$ with concentration $y$, which occurs with the same stoichiometry on both sides of each reaction in which it participates. Then $(\mathcal{R}', \mathcal{K}')$ admits an NPPO (resp., SPPO). 
\end{thm}

\begin{pf}
Fix the rate function $v \in \mathcal{K}$ such that $\mathcal{R}$ has an NPPO (resp., SPPO) $\mathcal{O}$. Fix $w \in \mathcal{K}'$ such that $w(x,1) = v(x)$, possible by assumption. With this rate function, the evolution of $\mathcal{R}'$ is governed by 
\begin{equation}
\label{eqtriv}
\left(\begin{array}{c}\dot x\\ \dot y\end{array}\right) = \left(\begin{array}{c}\Gamma\\0\end{array}\right)w(x,y).
\end{equation}
(\ref{eqtriv}) leaves $\mathbb{R}^n_{\gg 0} \times \{y\}$ locally invariant for each $y >0$. Since $w(x,1) = v(x)$, $\mathcal{O}':= \mathcal{O} \times \{1\}$ is a periodic orbit of $\mathcal{R}'$. Let $S' = S \times \{0\}$ so that $S'_\mathcal{O'}:=(\mathcal{O}' + S') \cap (\mathbb{R}^{n}_{\gg 0} \times \mathbb{R}_{> 0}) = S_\mathcal{O} \times \{1\}$ is the positive stoichiometry class of $\mathcal{O}'$ for $\mathcal{R}'$. Define $h \colon \mathbb{R}^r \to x_0+S$ by $h(z) = x_0+\Gamma_0 z$. $h$ is an affine bijection between $W:=h^{-1}(S_\mathcal{O}) \subseteq \mathbb{R}^r$ and $S_\mathcal{O}$, and defines local coordinates on $S_\mathcal{O}$ which evolve according to
\[
\dot z = Qv(x_0+\Gamma_0 z)\,.
\]
By definition, as $\mathcal{O}$ is an NPPO (resp., SPPO), $h^{-1}(\mathcal{O})$ is nondegenerate (resp., linearly stable). Now define $h' \colon \mathbb{R}^r \to (x_0+S) \times \{1\}$ by $h'(z) = (x_0+\Gamma_0 z, 1)$ and note that $h'$ is an affine bijection between $W$ and $S'_\mathcal{O'}$. Moreover $h'$ gives rise to precisely the same evolution in local coordinates (since $w(x_0+\Gamma_0 z, 1) = v(x_0 + \Gamma_0 z)$) and $h'^{-1}(\mathcal{O}') = h^{-1}(\mathcal{O})$. Thus, by definition, $\mathcal{O}'$ is an NPPO (resp., SPPO) of $\mathcal{R}'$. 
\hfill$\square$
\end{pf}

\begin{thm}[Adding a new species with inflow and outflow]
\label{thmnewwithopen}
Let $(\mathcal{R}, \mathcal{K})$ be a CRN with $C^2$ kinetics admitting an NPPO (resp., SPPO). Let $(\mathcal{R}', \mathcal{K}')$ be a species-reaction-extension of $(\mathcal{R}, \mathcal{K})$ created by
\begin{enumerate}[align=left, leftmargin=*]
\item[(i)] adding into the reactions of $\mathcal{R}$ a new species $Y$ with arbitrary stoichiometries; and
\item[(ii)] adding the new reaction $0 \rightleftharpoons Y$ with $C^2$ kinetics belonging to a scaling invariant subset of positive general kinetics.
\end{enumerate}
Then $(\mathcal{R}', \mathcal{K}')$ admits an NPPO (resp., SPPO).
\end{thm}

\begin{pf}
Fix the rate function $v \in \mathcal{K}$ such that $\mathcal{R}$ has an NPPO (resp., SPPO) $\mathcal{O}$. As in the proof of Theorem~\ref{thmnewdepreac} define $h \colon \mathbb{R}^r \to x_0+S$ by $h(z) = x_0+\Gamma_0 z$ and note that $h$ is an affine bijection between the open set $h^{-1}(S_\mathcal{O}) \subseteq \mathbb{R}^r$ and $S_\mathcal{O}$. $h$ defines local coordinates on $S_\mathcal{O}$ via $x = h(z)$, and $z$ evolves according to 
\begin{equation}
\label{eqbasic}
\dot z = Qv(x_0+\Gamma_0 z)\,.
\end{equation}
(\ref{eqbasic}) has a hyperbolic (resp., linearly stable) periodic orbit $\mathcal{O}' = h^{-1}(\mathcal{O})$. The assumptions on the kinetics mean that: 
\begin{enumerate}[align=left,leftmargin=*]
\item The new rate function $w(x,y)$ of the existing reactions can be chosen to satisfy $w(x,1) = v(x)$.
\item There exists a $C^2$ function $f\colon \mathbb{R}_{>0} \to \mathbb{R}_{>0}$ satisfying $f(1)=1$ and $f'(y)>0$ for all $y > 0$ and such that we may choose the rate of $0 \rightleftharpoons Y$ to be $\frac{1}{\epsilon}(1-f(y))$ where $\epsilon > 0$ is a parameter to be controlled. 
\end{enumerate}
With these choices, $\mathcal{R}'$ gives rise to the following singularly perturbed system:
\begin{equation}
\specialnumber{${}_\epsilon$}\label{eqSP}
\begin{array}{rcl}\dot x & = & \Gamma w(x,y)\\\epsilon\dot y & = & \epsilon s w(x,y) + (1-f(y)).\end{array}
\end{equation}
Here $s_i$ is the net change in the stoichometry of $Y$ in the $i$th reaction of $\mathcal{R}'$, and $s:=(s_1, \ldots, s_m)^{\mathrm{t}}$. For any fixed $\epsilon > 0$, rescaling time in the ``slow time system'' \specialeqref{eqSP}{${}_\epsilon$} gives the ``fast time system'':
\begin{equation}
\specialnumber{${}_\epsilon$}\label{eqSPa}
\begin{array}{rcl}\dot x & = & \epsilon\Gamma w(x,y)\\\dot y & = & \epsilon s w(x,y) + (1-f(y)).\end{array}
\end{equation}
Define $h' \colon \mathbb{R}^r\times \mathbb{R} \to (x_0+S) \times \mathbb{R}$ by $h'(z,y)  := (h(z), y)= (x_0+\Gamma_0 z, y)$. Note that $h'$ is an affine bijection between $h^{-1}(S_\mathcal{O}) \times \mathbb{R}_{> 0}$ and $S_\mathcal{O} \times \mathbb{R}_{> 0}$ and defines local coordinates on $S_\mathcal{O} \times \mathbb{R}_{> 0}$ via $(x,y) = (h(z),y)$. In $(z,y)$ coordinates the slow time system \specialeqref{eqSP}{${}_\epsilon$} becomes:
\begin{equation}
\specialnumber{${}_\epsilon$}\label{eqSPred}
\begin{array}{rcl}\dot z&=&Qw(x_0+\Gamma_0 z,y)\\\epsilon\dot y&=&\epsilon sw(x_0+\Gamma_0 z,y) + (1-f(y)).\end{array}
\end{equation}
while the fast time system \specialeqref{eqSPa}{${}_\epsilon$} becomes:
\begin{equation}
\specialnumber{${}_\epsilon$}\label{eqSPreda}
\begin{array}{rcl}\dot z&=&\epsilon Q w(x_0+\Gamma_0 z,y)\\\dot y&=&\epsilon s w(x_0+\Gamma_0 z,y) + (1-f(y)).\end{array}
\end{equation}
\specialeqref{eqSPred}{${}_0$} is the decoupled differential-algebraic system
\[
\dot z = Q v(x_0+\Gamma_0 z) =: H(z), \quad y=1\,,
\]
(as $f(y) = 1$ if and only if $y=1$, and $w(x_0+\Gamma_0 z, 1) = v(x_0+\Gamma_0 z)$). We observe that
\begin{enumerate}[align=left,leftmargin=*]
\item The vector field $H(z)$ has a hyperbolic (resp., linearly stable) periodic orbit $\mathcal{O}'$ by assumption, and hence \specialeqref{eqSPred}{${}_0$} has a periodic orbit $\overline{\mathcal{O}}:=\mathcal{O}'\times\{1\}$.
\item $y=1$ is a linearly stable equilibrium of $\dot y = 1-f(y)$ or, equivalently, the Jacobian matrix of \specialeqref{eqSPreda}{${}_\epsilon$} evaluated at $y=1, \epsilon = 0$, namely,
\[
\left(\begin{array}{cc}0&0\\0&-f'(1)\end{array}\right)\,,
\]
has a single nontrivial eigenvalue $-f'(1) < 0$. 
\end{enumerate}
By Theorems~13.1~and~13.2 in \cite{Fenichel79}, observations~(1)~and~(2) together tell us that there exists $\epsilon_0 > 0$ s.t. that for $\epsilon \in (0, \epsilon_0)$, \specialeqref{eqSPred}{${}_\epsilon$} has a hyperbolic (resp., linearly stable) periodic orbit $\overline{\mathcal{O}}_\epsilon$ close to $\overline{\mathcal{O}}$. Thus, for $\epsilon \in (0, \epsilon_0)$, $\mathcal{R}'$ has an NPPO (resp., SPPO) $\mathcal{O}_\epsilon := h'(\overline{\mathcal{O}}_\epsilon)$. 
\hfill$\square$
\end{pf}

\begin{remark}[Geometrical interpretation of the proof of Theorem~\ref{thmnewwithopen}]
The differential algebraic system \specialeqref{eqSPred}{${}_0$} defines a local flow on the $r$-dimensional (smooth) manifold $\mathcal{Y} := h^{-1}(\mathcal{S}_{\mathcal{O}}) \times \{1\}$ which includes the periodic orbit $\overline{\mathcal{O}}$. Let $\mathcal{Y}_0$ be some compact subset of $\mathcal{Y}$ containing $\overline{\mathcal{O}}$. The theory developed by Fenichel \cite{Fenichel79} shows (roughly, and omitting a myriad of technical details) that for sufficiently small $\epsilon > 0$ \specialeqref{eqSPred}{${}_\epsilon$} has an $r$-dimensional locally invariant manifold $\mathcal{Y}_\epsilon$ close to $\mathcal{Y}_0$. The vector field of \specialeqref{eqSPred}{${}_\epsilon$} restricted to $\mathcal{Y}_\epsilon$ is $\epsilon$-close to that of \specialeqref{eqSPred}{${}_0$} on $\mathcal{Y}_0$ and consequently, by regular perturbation theory, for sufficiently small $\epsilon \neq 0$, \specialeqref{eqSPred}{${}_\epsilon$} has a periodic orbit $\overline{\mathcal{O}}_\epsilon$ on $\mathcal{Y}_\epsilon$ close to $\overline{\mathcal{O}}$. The technical assumption that all vector fields involved are $C^2$ is to ensure that the family of vector fields on $\mathcal{Y}_\epsilon$ is $C^1$, allowing use of regular perturbation theory. The $r-1$ nontrivial Floquet multipliers of $\overline{\mathcal{O}}_\epsilon$ relative to $\mathcal{Y}_\epsilon$ are close to those of $\overline{\mathcal{O}}$ relative to $\mathcal{Y}$ which, by assumption, are disjoint from (resp., inside) the unit circle. Meanwhile, the single Floquet multiplier of $\overline{\mathcal{O}}_\epsilon$ transverse to $\mathcal{Y}_\epsilon$ lies inside the unit circle as a consequence of the fact that $-f'(1) < 0$.  
\end{remark}

\begin{remark}[Kinetic assumptions in Theorem~\ref{thmnewwithopen}]
\label{remnewwithopen}
The added flow reaction $0 \rightleftharpoons Y$ in Theorem~\ref{thmnewwithopen} may have, for example, mass action kinetics, positive general kinetics, physical power-law kinetics, or any fixed physical power-law kinetics (these all define scaling invariant subsets of positive general kinetics). 
\end{remark}

Theorem~\ref{thmnewwithopen}, combined with Theorem~\ref{thmnewdepreac} allows us to deduce an important corollary:
\begin{prop}[Inheritance in fully open species-reaction extensions]
\label{coropeninduced}
Let $(\mathcal{R}, \mathcal{K})$ be a fully open CRN with $C^2$ kinetics admitting an NPPO (resp., SPPO). Let $(\mathcal{R}', \mathcal{K}')$ be a fully open CRN with kinetics, which is a species-reaction extension of $(\mathcal{R}, \mathcal{K})$ (Definition~\ref{defspecreacext}), and such that for each new reaction $R$ in $\mathcal{R}'$, $\mathcal{K}'^{(R)}$ is $C^2$, and belongs to a scaling invariant subset of positive general kinetics. Then $(\mathcal{R}', \mathcal{K}')$ admits an NPPO (resp., SPPO).
\end{prop}
\begin{pf}
Let $\mathcal{R}$ have $n_1$ species and $m_1$ non-flow reactions (i.e., reactions not of the form $0 \rightarrow X_i$ or $X_i \rightarrow 0$), and $\mathcal{R}'$ have $n_2$ species and $m_2$ non-flow reactions. We can construct $(\mathcal{R}', \mathcal{K}')$ from $(\mathcal{R}, \mathcal{K})$ via a sequence of steps as follows:
\begin{enumerate}
\item[(i)] Beginning with $\mathcal{R}$, for each absent species $X_j$ (if any) we add the species to all existing reactions and add $0 \rightleftharpoons X_j$. The kinetic assumptions ensure that this corresponds to $n_2-n_1$ applications of Theorem~\ref{thmnewwithopen}, one for each absent species. Note that, as $\mathcal{R}$ is fully open, the new CRN created at each stage is fully open.
\item[(ii)] We add each remaining absent reaction (if any). The kinetic assumptions ensure that this corresponds to $m_2-m_1$ applications of Theorem~\ref{thmnewdepreac}, one for each reaction added. Theorem~\ref{thmnewdepreac} applies because a fully open CRN has stoichiometric subspace which is the whole state space, and hence any added reaction is a dependent reaction. 
\end{enumerate}
We can see the above procedure as constructing a sequence of intermediate (fully-open) CRNs with kinetics, beginning with $(\mathcal{R}, \mathcal{K})$ and terminating with $(\mathcal{R}', \mathcal{K}')$:
\[
(\mathcal{R}, \mathcal{K}) = (\mathcal{R}_0, \mathcal{K}_0)\underbrace{\strut \quad \rightarrow \quad \cdots \quad \rightarrow \quad }_{\mathclap{\mbox{add in species and flows (Thm.~\ref{thmnewwithopen})}}}(\mathcal{R}_{n_2-n_1}, \mathcal{K}_{n_2-n_1})\underbrace{\strut \quad \rightarrow \quad \cdots \quad \rightarrow \quad }_{\mathclap{\mbox{add in reactions (Thm.~\ref{thmnewdepreac})}}}(\mathcal{R}_p, \mathcal{K}_p) = (\mathcal{R}', \mathcal{K}')\,.
\]
(Here $p = n_2+m_2-n_1-m_1$.) If $(\mathcal{R}, \mathcal{K})$ admits an NPPO (resp., SPPO), then each step of the above procedure preserves this property, and consequently $(\mathcal{R}', \mathcal{K}')$ admits an NPPO (resp., SPPO).
\hfill$\square$
\end{pf}

\begin{remark}[Kinetic assumptions in Proposition~\ref{coropeninduced}, proof of Proposition~\ref{propMAfo}]
\label{remMAfo}
The somewhat unwieldy kinetic assumptions in Proposition~\ref{coropeninduced} are in order to maximise generality. They are satisfied if $\mathcal{R} \leq \mathcal{R}'$ and, for example, 
\begin{enumerate}[align=left,leftmargin=*]
\item Both $\mathcal{R}$ and $\mathcal{R}'$ have mass action kinetics. 
\item Both $\mathcal{R}$ and $\mathcal{R}'$ have physical power-law kinetics. 
\item $\mathcal{R}$ has any fixed power-law kinetics and $\mathcal{R}'$ has any power-law kinetics derived from that of $\mathcal{R}$ (see Definition~\ref{derived}).
\item Both $\mathcal{R}$ and $\mathcal{R}'$ have $C^2$ positive general kinetics. 
\end{enumerate}
Thus, in particular Proposition~\ref{propMAfo} follows immediately from Proposition~\ref{coropeninduced}. Explorations in Section~\ref{secnum} are carried out using Proposition~\ref{coropeninduced} with $\mathcal{R}$ and $\mathcal{R}'$ both given mass action kinetics or both given physical power-law kinetics. 
\end{remark}

In the light of Proposition~\ref{coropeninduced}, and adapting the terminology of \cite{joshishiu}, the following definitions make sense.
\begin{def1}[Atoms of oscillation, atoms of stable oscillation]
\label{atomosci}
A fully open mass action CRN which admits an NPPO (resp., SPPO), and which is minimal with respect to the induced subnetwork ordering amongst fully open mass action CRNs admitting NPPOs (resp., SPPOs), is referred to as a {\em fully open mass action atom of oscillation} (resp., {\em stable oscillation}). Atoms with respect to other classes of kinetics, such as physical power-law kinetics, are similarly defined. 
\end{def1}

Observe that Definition~\ref{atomosci} is restricted to {\em fully open} CRNs as the presence of an oscillatory induced subnetwork in a general CRN does not necessarily imply oscillation; an example is provided in the concluding section (Example~\ref{exinherit}). Note also that, as in the case of multistationarity \cite{banajipanteaMPNE}, a fully open mass action atom of oscillation with respect to the induced subnetwork ordering may not be minimal with respect to other, better partial orders. Note finally that a fully open mass action atom of oscillation may include an induced subnetwork admitting an NPPO but which is not fully open; thus if we do not restrict attention to fully open CRNs, fully open atoms of oscillation may not be minimal oscillatory CRNs even with respect to the induced subnetwork ordering.

\section{The occurrence of stable oscillation in small, fully open, CRNs}
\label{secnum}

A fully open CRN is taken to be ``small'' if it has few species, few non-flow reactions, and is at most bimolecular, namely the total stoichiometry of all species on each side of every reaction is no more than two. The goal of this section is to provide some lower bounds on the frequency with which small fully open CRNs admit SPPOs under the assumptions of (i) mass action kinetics and (ii) physical power-law kinetics. This is done via a mixture of basic analysis, numerical simulation, and application of the inheritance result in Proposition~\ref{coropeninduced}. 

Define a $(k,l)$ CRN to be a fully open, at most bimolecular, CRN with $k\geq 1$ species and $l\geq 0$ irreversible non-flow reactions. It is easy to see that $(1,l)$ and $(k,0)$ CRNs can admit no nontrivial periodic orbits for any reasonable kinetics: if $k=1$ then regardless of the kinetics (\ref{genCRN}) is a one dimensional autonomous system which forbids nontrivial oscillation; if $l=0$ then, with positive general kinetics, (\ref{genCRN}) is a decoupled system of $k$ autonomous univariate ODEs which again forbids nontrivial oscillation. 

We proceed as follows. We first treat the smallest nontrivial case, namely $(k,l)=(2,1)$, which is simple enough to be fully analysed using fairly basic ideas from dynamical systems. The results of this analysis are summarised in Propostion~\ref{prop21}. We then proceed as follows, ensuring that $(k,l)$ CRNs are treated after $(k-1,l)$ and $(k,l-1)$ CRNs, and treating the cases of mass action, and of physical power-law kinetics separately. 

\begin{enumerate}[align=left,leftmargin=*]
\item Whenever an SPPO is found in a $(k,l)$ CRN $\mathcal{R}$, we use the powerful and widely available graph-isomorphism software NAUTY \cite{nauty} to identify all $(k+1,l)$ CRNs and $(k,l+1)$ CRNs $\geq \mathcal{R}$ (i.e., which include $\mathcal{R}$ as an induced subnetwork). Proposition~\ref{coropeninduced} then tells us that these must admit SPPOs.
\item We search numerically for SPPOs in $(k,l)$ CRNs, limiting the search to those CRNs not already found to inherit SPPOs via step (1) or believed to forbid oscillation by Conjecture~\ref{conjJac} below.
\end{enumerate}

Via this process we obtain a lower bound on the occurrence of stable oscillation in small fully open CRNs. The methodological details are in \ref{appmethod}.

\begin{prop}
\label{prop21}
There are 14 non-isomorphic $(2,1)$ CRNs. These are, upto isomorphism, the fully open extensions of:
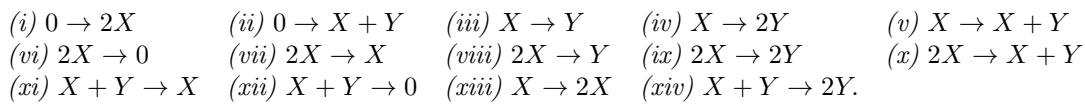
\[
\begin{array}{lllll}
\mbox{(i)}\,\,0\rightarrow 2X & \mbox{(ii)}\,\,0\rightarrow X+Y & \mbox{(iii)}\,\,X\rightarrow Y & \mbox{(iv)}\,\,X\rightarrow 2Y & \mbox{(v)}\,\,X\rightarrow X+Y\\ \mbox{(vi)}\,\,2X\rightarrow 0 & \mbox{(vii)}\,\,2X\rightarrow X & \mbox{(viii)}\,\,2X\rightarrow Y & \mbox{(ix)}\,\,2X\rightarrow 2Y & \mbox{(x)}\,\,2X\rightarrow X+Y\\ \mbox{(xi)}\,\,X+Y\rightarrow X & \mbox{(xii)}\,\,X+Y\rightarrow 0 & \mbox{(xiii)}\,\,X\rightarrow 2X & \mbox{(xiv)}\,\,X+Y\rightarrow 2Y.&
\end{array}
\]
Let $\mathcal{R}_{(k)}$ refer to the fully open extension of reaction (k), namely the CRN consisting of reaction (k) along with $X\rightleftharpoons 0 \rightleftharpoons Y$.
\begin{enumerate}[align=left,leftmargin=*]
\item With mass action kinetics $\mathcal{R}_{(i)}$ to $\mathcal{R}_{(xiv)}$ forbid oscillation. All but $\mathcal{R}_{(xiii)}$ have a unique equilibrium which is locally asymptotically stable and attracts all of $\mathbb{R}^2_{\geq 0}$. $\mathcal{R}_{(xiii)}$ either has a unique locally asymptotically stable equilibrium which attracts all of $\mathbb{R}^2_{\geq 0}$, or all orbits are unbounded.
\item With positive general kinetics $\mathcal{R}_{(i)}$ to $\mathcal{R}_{(xiii)}$ forbid oscillation.
\item With physical power-law kinetics or general kinetics $\mathcal{R}_{(xiv)}$ admits an SPPO.
\end{enumerate}
\end{prop}

The proof of Proposition~\ref{prop21} is fairly straightforward, but somewhat lengthy, and is in \ref{pf21}. In order to proceed more efficiently, we make the following conjecture.
\begin{conj}
\label{conjJac}
Let $(\mathcal{R}, \mathcal{K})$ be a fully open CRN with kinetics, where $\mathcal{K}$ is any scaling invariant subset of positive general kinetics (for example, $\mathcal{K}$ may be given by mass action kinetics or physical power-law kinetics). Suppose that $\Gamma$ is the stoichiometric matrix of $\mathcal{R}$ so that $\mathcal{R}$ gives rise to the family of ODEs on $\mathbb{R}^n_{\gg 0}$
\[
\dot x = \Gamma v(x), \quad v \in \mathcal{K}\,.
\]
If, for all $x \in \mathbb{R}^n_{\gg 0}$ and all $v \in \mathcal{K}$, the Jacobian matrix $\Gamma Dv(x)$ has no purely imaginary eigenvalues, then $\mathcal{R}$ does not admit a positive periodic orbit. 
\end{conj}

A theoretical justification for Conjecture~\ref{conjJac} is not attempted here, but it is not hard to believe the rather stronger claim that such families of CRNs admit oscillation if and only if they admit Hopf bifurcation (a similar conjecture is made in Section~2.2 of \cite{eiswirth91}). If Conjecture~\ref{conjJac} holds, it is possible to rule out oscillation by examining, with the help of computer algebra, certain polynomials associated with $\Gamma Dv(x)$ whose positivity is sufficient to forbid purely imaginary eigenvalues. This process, which will be described in forthcoming work, is computationally much less expensive than simulating the differential equations with tens of thousands of parameter choices. No counterexamples to Conjecture~\ref{conjJac} were found during a large number of numerical simulations. Note also that as the claims such as those drawn from the data in Table~\ref{tabdat} concern {\em lower bounds} on the frequency of oscillation in CRNs, they are not invalidated if Conjecture~\ref{conjJac} is false.

\begin{table}[h]
\begin{center}
\begin{tikzpicture}[scale=0.5]

\fill[color=black!20] (9,3) -- (15,3) -- (15,6) -- (9,6) -- cycle;

\draw [-, line width=0.04cm] (2,11) -- (23,11);
\draw [-, line width=0.04cm] (2,10) -- (23,10);
\draw [-, line width=0.04cm] (0,9) -- (23,9);
\draw [-, line width=0.04cm] (1,6) -- (23,6);
\draw [-, line width=0.04cm] (1,3) -- (23,3);
\draw [-, line width=0.04cm] (0,0) -- (23,0);

\draw [-, line width=0.04cm] (0,0) -- (0,9);
\draw [-, line width=0.04cm] (1,0) -- (1,9);
\draw [-, line width=0.04cm] (2,0) -- (2,11);
\draw [-, line width=0.04cm] (5,0) -- (5,10);
\draw [-, line width=0.04cm] (9,0) -- (9,10);
\draw [-, line width=0.04cm] (15,0) -- (15,10);
\draw [-, line width=0.04cm] (23,0) -- (23,11);

\draw [-, line width=0.01cm] (2,1) -- (23,1);
\draw [-, line width=0.01cm] (2,2) -- (23,2);
\draw [-, line width=0.01cm] (2,4) -- (23,4);
\draw [-, line width=0.01cm] (2,5) -- (23,5);
\draw [-, line width=0.01cm] (2,7) -- (23,7);
\draw [-, line width=0.01cm] (2,8) -- (23,8);

\draw [-, line width=0.01cm] (3.5,0) -- (3.5,2);
\draw [-, line width=0.01cm] (3.5,3) -- (3.5,5);
\draw [-, line width=0.01cm] (3.5,6) -- (3.5,8);

\draw [-, line width=0.01cm] (7,0) -- (7,2);
\draw [-, line width=0.01cm] (7,3) -- (7,5);
\draw [-, line width=0.01cm] (7,6) -- (7,8);

\draw [-, line width=0.01cm] (12,0) -- (12,2);
\draw [-, line width=0.01cm] (12,3) -- (12,5);
\draw [-, line width=0.01cm] (12,6) -- (12,8);

\draw [-, line width=0.01cm] (19,0) -- (19,2);
\draw [-, line width=0.01cm] (19,3) -- (19,5);
\draw [-, line width=0.01cm] (19,6) -- (19,8);

\node at (12.5,10.5) {number of non-flow reactions $l$};
\node[rotate=90] at (0.5,4.5) {number of species $k$};

\node at (3.5,9.5) {$1$};
\node at (7,9.5) {$2$};
\node at (12,9.5) {$3$};
\node at (19,9.5) {$4$};

\node at (1.5,7.5) {$2$};
\node at (1.5,4.5) {$3$};
\node at (1.5,1.5) {$4$};

\node at (3.5,8.5) {14};
\node at (7,8.5) {169};
\node at (12,8.5) {1,312};
\node at (19,8.5) {7,514};

\node at (3.5,5.5) {19};
\node at (7,5.5) {622};
\node at (12,5.5) {16,135};
\node at (19,5.5) {322,854};

\node at (3.5,2.5) {20};
\node at (7,2.5) {1,059};
\node at (12,2.5) {59,379};
\node at (19,2.5) {2,840,062};
\node at (2.75,7.5) {0};
\node at (4.25,7.5) {0};
\node at (6,7.5) {0};
\node at (8,7.5) {0};
\node at (10.5,7.5) {0};
\node at (13.5,7.5) {0};
\node at (17,7.5) {0};
\node at (21,7.5) {0};
\node at (2.75,6.5) {1};
\node at (4.25,6.5) {0};
\node at (6,6.5) {25};
\node at (8,6.5) {25};
\node at (10.5,6.5) {293};
\node at (13.5,6.5) {289};
\node at (17,6.5) {2,257};
\node at (21,6.5) {2,246};
\node at (2.75,4.5) {0};
\node at (4.25,4.5) {0};
\node at (6,4.5) {5};
\node at (8,4.5) {0};
\node at (10.5,4.5) {444};
\node at (13.5,4.5) {401};
\node at (17,4.5) {$\geq$ 18,859};
\node at (21,4.5) {18,859};
\node at (2.75,3.5) {1};
\node at (4.25,3.5) {1};
\node at (6,3.5) {94};
\node at (8,3.5) {82};
\node at (10.5,3.5) {4,268};
\node at (13.5,3.5) {4,080};
\node at (17,3.5) {$\geq$ 123,990};
\node at (21,3.5) {123,990};
\node at (2.75,1.5) {0};
\node at (4.25,1.5) {0};
\node at (6,1.5) {8};
\node at (8,1.5) {8};
\node at (10.5,1.5) {$\geq$ 1,657};
\node at (13.5,1.5) {1,657};
\node at (17,1.5) {$\geq$ 166,676};
\node at (21,1.5) {166,676};
\node at (2.75,0.5) {1};
\node at (4.25,0.5) {1};
\node at (6,0.5) {140};
\node at (8,0.5) {139};
\node at (10.5,0.5) {$\geq$ 14,373};
\node at (13.5,0.5) {14,373};
\node at (17,0.5) {$\geq$ 1,038,785};
\node at (21,0.5) {1,038,785};

\end{tikzpicture}
\end{center}

\caption{\label{tabdat} The table shows (i) the total number of nonisomorphic $(k,l)$ CRNs for $k = 2,\ldots, 4$ and $l=1, \ldots, 4$, (ii) lower bounds on the number of $(k,l)$ CRNs admitting SPPOs under the assumptions of mass action kinetics and physical power-law kinetics, and (iii) lower bounds on how many of these admit SPPOs as a consequence of the inheritance results in this paper. Each block of five cells corresponding to a particular pair of $(k,l)$ contains the total number of nonisomorphic $(k,l)$ CRNs (top row); the number shown to admit SPPOs with mass action kinetics followed by the number of these which follow as a consequence of inheritance results (middle row); and the number shown to admit SPPOs with physical power-law kinetics followed by the number of these which follow as a consequence of inheritance results (bottom row). For example, the data in the highlighted block tells us that there are 16,135 nonisomorphic $(3,3)$ CRNs. Of these, at least 444 (about $3\%$) admit SPPOs with mass action kinetics: 401 (about $90\%$) by inheritance, namely because they include as an induced subnetwork either a $(3,2)$ CRN or a $(2,3)$ CRN which admits an SPPO, with the remainder found in numerical simulations. Similarly, at least 4,264 (about $26\%$) of the $(3,3)$ CRNs admit SPPOs with physical power-law kinetics: 4,072 (about $95\%$) by inheritance, with the remainder being found in numerical simulations. For $k+l \geq 7$, only the inheritance data is presented namely, no numerical search was carried out to find CRNs admitting SPPOs not predicted by the inheritance results. A ``$\geq$'' is inserted in order to highlight this. The lists of CRNs from which the data is drawn are at \protect\url{https://reaction-networks.net/networks/osci.html}.
}
\end{table}

The results of simulations and analysis for $k= 2,\ldots, 4$ and $l=1, \ldots, 4$ are summarised in Table~\ref{tabdat}. The table suggests, assuming that Conjecture~\ref{conjJac} is true, and that large numbers of oscillatory CRNs were not missed by the numerical simulations, that the great majority of CRNs admitting stable oscillation do so as a consequence of inheritance (this becomes even more evident as we increase the number of reactions in the CRNs). As a particular example, the motif 
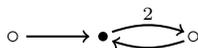
\begin{center}
\begin{tikzpicture}[scale=1.2]
\fill[color=black] (3,0) circle (1.5pt);
\draw (2,0) circle (1.5pt);
\draw (4,0) circle (1.5pt);


\draw [->, thick] (2.15,0) -- (2.85,0);
\draw [->, thick] (3.1,0.05) .. controls (3.4,0.15) and (3.6,0.15) .. (3.9,0.05);
\draw [->, thick] (3.9,-0.05) .. controls (3.6,-0.15) and (3.4,-0.15) .. (3.1,-0.05);

\node at (3.5,0.25) {$\scriptstyle{2}$};
\end{tikzpicture}
\end{center}
corresponding to the single reaction $X+Y \rightarrow 2Y$ occurs in $22\%$
of all the CRNs in Table~\ref{tabdat}, which consequently admit SPPOs with physical power-law kinetics by Propositions~\ref{prop21}~and~\ref{coropeninduced}. A total of about $33\%$ of the CRNs in Table~\ref{tabdat} were found to admit SPPOs with physical power-law kinetics
and thus this single motif is responsible for about two thirds of the oscillation found under the assumption of physical power-law kinetics. Additional investigation revealed that this motif occurs in a total of about $2.52\times 10^7$ 
($75\%$) of all 
$3.36 \times 10^7$ $(2,l)$ CRNs ($l$ ranges from $1$ to $26$ by the counting arguments in \ref{appmethod}). Thus identifying small atoms of oscillation is worthwhile from a practical viewpoint, as these appear to be the source of most oscillation in CRNs. 

Table~\ref{tabdat} also highlights the importance of kinetics, and in particular how much more frequently stable oscillation occurs in small CRNs with physical power-law kinetics as compared to those with mass action kinetics. Presumably the linear or quadratic nature of at most bimolecular mass action systems significantly restricts the allowed dynamics in many cases. 

As in the case of physical power-law kinetics, small oscillatory motifs account for most of the oscillation in the table found in mass action CRNs. For example, at least one of the five $(3,2)$ (presumed) mass action atoms of stable oscillation found in simulations occurs in about $5\%$
of all the CRNs in Table~\ref{tabdat}, which consequently admit SPPOs with mass action kinetics; this accounts for almost $90\%$ of the oscillation in mass action CRNs detailed in Table~\ref{tabdat}. While numerical investigations in \ref{appmethod} indicate that the lower bounds in Table~\ref{tabdat} can be improved with additional simulation, it remains true that inheritance results applied to a few small oscillatory motifs automatically give us large numbers of oscillatory CRNs.

Not visible in the table are relationships amongst the atoms of stable oscillation. For example, the five $(3,2)$ mass action atoms of stable oscillation are the fully open extensions of: (i) $X+Z\rightarrow 2Y \rightarrow Y+Z$, (ii) $X+Z \rightarrow 2Y$, $Y+Z \rightarrow 2Z$, (iii) $X+Z \rightarrow Y, Y+Z \rightarrow 2Z$, (iv) $X+Z \rightarrow 2Y \rightarrow 2Z$; (v) $X+Z \rightarrow 0, Y+Z \rightarrow 2Z$. These correspond to the following motifs:
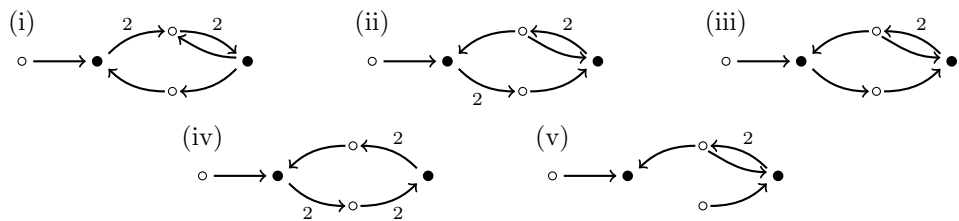
\begin{center}
\begin{tikzpicture}[scale=1]
\node at (1,0.5) {(i)};
\draw[color=black] (3,0.4) circle (1.5pt);
\draw[color=black] (3,-0.4) circle (1.5pt);
\draw[color=black] (1,0) circle (1.5pt);
\fill[color=black] (2,0) circle (2pt);
\fill[color=black] (4,0) circle (2pt);

\draw [->, thick] (1.15,0) --(1.85,0);
\draw [->, thick] (2.15,0.1) .. controls (2.3,0.25) and (2.5,0.4) .. (2.9,0.4);
\draw [<-, thick] (3.85,0.1) .. controls (3.7,0.25) and (3.5,0.4) .. (3.1,0.4);
\draw [<-, thick] (2.15,-0.1) .. controls (2.3,-0.25) and (2.5,-0.4) .. (2.9,-0.4);
\draw [->, thick] (3.85,-0.1) .. controls (3.7,-0.25) and (3.5,-0.4) .. (3.1,-0.4);
\draw[->, thick] (3.85,0.05) .. controls (3.5, 0.05) and (3.3, 0.2) .. (3.07,0.33);

\node at (3.6,0.5) {$\scriptstyle{2}$};
\node at (2.4,0.5) {$\scriptstyle{2}$};
\node at (2.4,-0.5) {$\textcolor{white}{\scriptstyle{1}}$};

\end{tikzpicture}
\hspace{1cm}
\begin{tikzpicture}[scale=1]

\node at (1,0.5) {(ii)};
\draw[color=black] (3,0.4) circle (1.5pt);
\draw[color=black] (3,-0.4) circle (1.5pt);
\draw[color=black] (1,0) circle (1.5pt);
\fill[color=black] (2,0) circle (2pt);
\fill[color=black] (4,0) circle (2pt);

\draw [->, thick] (1.15,0) --(1.85,0);
\draw [<-, thick] (2.15,0.1) .. controls (2.3,0.25) and (2.5,0.4) .. (2.9,0.4);
\draw [->, thick] (3.85,0.1) .. controls (3.7,0.25) and (3.5,0.4) .. (3.1,0.4);
\draw[<-, thick] (3.85,0.05) .. controls (3.5, 0.05) and (3.3, 0.2) .. (3.07,0.33);
\draw [->, thick] (2.15,-0.1) .. controls (2.3,-0.25) and (2.5,-0.4) .. (2.9,-0.4);
\draw [<-, thick] (3.85,-0.1) .. controls (3.7,-0.25) and (3.5,-0.4) .. (3.1,-0.4);

\node at (3.6,0.5) {$\scriptstyle{2}$};
\node at (2.4,-0.5) {$\scriptstyle{2}$};
\end{tikzpicture}
\hspace{1cm}
\begin{tikzpicture}[scale=1]

\node at (1,0.5) {(iii)};
\draw[color=black] (3,0.4) circle (1.5pt);
\draw[color=black] (3,-0.4) circle (1.5pt);
\draw[color=black] (1,0) circle (1.5pt);
\fill[color=black] (2,0) circle (2pt);
\fill[color=black] (4,0) circle (2pt);

\draw [->, thick] (1.15,0) --(1.85,0);
\draw [<-, thick] (2.15,0.1) .. controls (2.3,0.25) and (2.5,0.4) .. (2.9,0.4);
\draw [->, thick] (3.85,0.1) .. controls (3.7,0.25) and (3.5,0.4) .. (3.1,0.4);
\draw[<-, thick] (3.85,0.05) .. controls (3.5, 0.05) and (3.3, 0.2) .. (3.07,0.33);
\draw [->, thick] (2.15,-0.1) .. controls (2.3,-0.25) and (2.5,-0.4) .. (2.9,-0.4);
\draw [<-, thick] (3.85,-0.1) .. controls (3.7,-0.25) and (3.5,-0.4) .. (3.1,-0.4);

\node at (3.6,0.5) {$\scriptstyle{2}$};
\node at (2.4,-0.5) {$\textcolor{white}{\scriptstyle{2}}$};
\end{tikzpicture}
\begin{tikzpicture}[scale=1]

\node at (1,0.5) {(iv)};
\draw[color=black] (3,0.4) circle (1.5pt);
\draw[color=black] (3,-0.4) circle (1.5pt);
\draw[color=black] (1,0) circle (1.5pt);
\fill[color=black] (2,0) circle (2pt);
\fill[color=black] (4,0) circle (2pt);

\draw [->, thick] (1.15,0) --(1.85,0);
\draw [<-, thick] (2.15,0.1) .. controls (2.3,0.25) and (2.5,0.4) .. (2.9,0.4);
\draw [->, thick] (3.85,0.1) .. controls (3.7,0.25) and (3.5,0.4) .. (3.1,0.4);
\draw [->, thick] (2.15,-0.1) .. controls (2.3,-0.25) and (2.5,-0.4) .. (2.9,-0.4);
\draw [<-, thick] (3.85,-0.1) .. controls (3.7,-0.25) and (3.5,-0.4) .. (3.1,-0.4);

\node at (3.6,0.5) {$\scriptstyle{2}$};
\node at (2.4,-0.5) {$\scriptstyle{2}$};
\node at (3.6,-0.5) {$\scriptstyle{2}$};
\end{tikzpicture}
\hspace{1cm}
\begin{tikzpicture}[scale=1]

\node at (1,0.5) {(v)};
\draw[color=black] (3,0.4) circle (1.5pt);
\draw[color=black] (3,-0.4) circle (1.5pt);
\draw[color=black] (1,0) circle (1.5pt);
\fill[color=black] (2,0) circle (2pt);
\fill[color=black] (4,0) circle (2pt);

\draw [->, thick] (1.15,0) --(1.85,0);
\draw [<-, thick] (2.15,0.1) .. controls (2.3,0.25) and (2.5,0.4) .. (2.9,0.4);
\draw [->, thick] (3.85,0.1) .. controls (3.7,0.25) and (3.5,0.4) .. (3.1,0.4);
\draw[<-, thick] (3.85,0.05) .. controls (3.5, 0.05) and (3.3, 0.2) .. (3.07,0.33);
\draw [<-, thick] (3.85,-0.1) .. controls (3.7,-0.25) and (3.5,-0.4) .. (3.1,-0.4);

\node at (3.6,0.5) {$\scriptstyle{2}$};
\node at (2.4,-0.5) {$\textcolor{white}{\scriptstyle{2}}$};
\end{tikzpicture}
\end{center}
Representing these motifs pictorially highlights the close relationships between them. Observe that there are various subnetwork relationships between the motifs. For example, (v) is a subnetwork of (iii), but not an {\em induced} subnetwork of (iii), and hence oscillation in the fully open extension of (iii) cannot be predicted from that in the fully open extension of (v) using the theorems in this paper. There remains the possibility that there exists an inheritance result rather different from those in this paper which predicts oscillation in the fully open extension of (iii) from that in the fully open extension of (v). More generally, it seems likely that there are interesting theorems to be discovered on sufficient conditions for stable oscillation in mass action CRNs which might explain something about the structures of oscillatory motifs.  

\section{Conclusions}

Armed with the results in this paper one can predict the occurrence of oscillation in CRNs from its occurrence in smaller CRNs. Our main conclusion is:
\begin{quote}
Any CRN built from an oscillatory CRN via a sequence of modifications of the kind described in Theorems~\ref{thmnewdepreac}~to~\ref{thmnewwithopen} is again oscillatory.
\end{quote}
Here ``oscillatory'' may be taken to mean either ``which admits an NPPO'' or ``which admits an SPPO'', and the conclusion is valid under mild assumptions on the kinetics and for general CRNs (not necessarily fully open). We emphasised the consequence that a fully open, mass action, CRN which includes a fully open oscillatory subnetwork is itself oscillatory, illustrating how certain motifs are associated with oscillation in fully open CRNs. It was mentioned, however, that this particular conclusion does not extend to CRNs which are not fully open: such a CRN may include an oscillatory subnetwork but fail to be oscillatory. The following is a typical example:

\begin{example}
\label{exinherit}
Consider the following CRNs $\mathcal{R}$, $\mathcal{R}'$ and $\mathcal{R}''$ which satisfy $\mathcal{R} \leq_{S} \mathcal{R}' \leq_{R} \mathcal{R}''$:
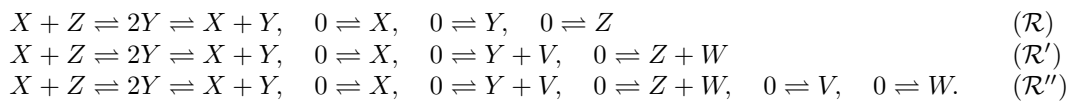
\[
\begin{array}{lcl}
X+Z\rightleftharpoons 2Y \rightleftharpoons X+Y, \quad 0 \rightleftharpoons X, \quad 0 \rightleftharpoons Y, \quad 0 \rightleftharpoons Z && (\mathcal{R})\\
X+Z\rightleftharpoons 2Y \rightleftharpoons X+Y, \quad 0 \rightleftharpoons X, \quad 0 \rightleftharpoons Y+V, \quad 0 \rightleftharpoons Z+W && (\mathcal{R}')\\
X+Z\rightleftharpoons 2Y \rightleftharpoons X+Y, \quad 0 \rightleftharpoons X, \quad 0 \rightleftharpoons Y+V, \quad 0 \rightleftharpoons Z+W, \quad 0 \rightleftharpoons V, \quad 0 \rightleftharpoons W. && (\mathcal{R}'')\\
\end{array}
\]
$\mathcal{R}$ admits an SPPO with mass action kinetics as it is just the fully open extension of motif (i) above, with the reverse of some reactions added (see Remark~\ref{newdeprev}). On the other hand $\mathcal{R}'$ is a weakly reversible, deficiency zero, network and, consequently, with mass action kinetics, forbids oscillation by the deficiency zero theorem \cite{feinberg}. Finally, by Theorem~\ref{thmnewwithopen} applied twice to $\mathcal{R}$, $\mathcal{R}''$ admits an SPPO with mass action kinetics. 
\end{example}

Example~\ref{exinherit} illustrates that predicting oscillation in CRNs is rather subtle: enlarging a CRN in natural ways can both destroy and create oscillation. Moreover, $\mathcal{R}'$ and $\mathcal{R}''$ involve the same set of species and have the same stoichiometric subspace (namely, all of $\mathbb{R}^5$); but adding the flow reactions $0 \rightleftharpoons V, \,0 \rightleftharpoons W$ to $\mathcal{R}'$ gives rise to oscillation. This corresponds to adding constant and linear terms to the differential equations describing the evolution of $\mathcal{R}'$ with mass action kinetics.

It is highly likely that further results of the kind presented in this paper hold: following Theorems~5~and~6 in \cite{banajipanteaMPNE} we expect that modifications such as including new reactions with new species, or inserting intermediate complexes involving new species into reactions should, with mild additional hypotheses, preserve the capacity for NPPOs or SPPOs. Some oscillatory CRNs, minimal w.r.t. to the modifications described in Theorems~\ref{thmnewdepreac}~to~\ref{thmnewwithopen} of this paper, may cease to be minimal under the improved partial order such results would bring.

There are also interesting questions on the connections between inheritance approaches as described here, and known sufficient conditions for oscillation such as those in \cite{eiswirth91, eiswirth96,gatermann, errami2015}. The families of chemical oscillators described in these papers can provide a starting point for application of the inheritance results here. It is also possible that some of the theory on families of chemical oscillators or algorithmic conditions for oscillation might suggest further inheritance results not described here. These possibilities remain to be explored.

The investigation of small, fully open, CRNs in Section~\ref{secnum} highlights two important points:
\begin{itemize}
\item identifying small oscillatory motifs is a worthwhile pursuit as it automatically implies oscillation in the large number of CRNs which ``inherit'' these motifs; and 
\item stable oscillation is much more common given larger classes of kinetics such as physical power-law kinetics as compared to mass action kinetics. 
\end{itemize}
Similar studies could also be carried out for general CRNs (not necessarily fully open), using Theorems~\ref{thmnewdepreac}~to~\ref{thmnewwithopen}. The difficulty of finding oscillation in mass action CRNs by numerical experiment is evidenced by additional data in \ref{appmethod}. This data suggests that often oscillation is confined to small parameter regions, and encourages the use of more systematic algorithmic approaches to the detection of oscillation such as those in \cite{errami2015}.

Finally, the ``enumerate and simulate'' methodology which provided the data in Section~\ref{secnum} and is described in more detail in \ref{appmethod} may also prove useful for studying the frequency of other behaviours such as chaos in CRNs \cite{pojman}. Some modification to the approach may be needed to explore sets of CRNs too large to be studied exhaustively. For example, there are more than $10^8$
nonisomorphic $(4,5)$ CRNs, and exploring the dynamics of such large numbers numerically becomes challenging; however, it should be possible either to restrict attention to certain interesting subsets of these CRNs, such as those which are weakly reversible for example, or to explore randomly chosen CRNs from such sets in order to draw some conclusions about how often various behaviours might occur. 

\section*{Acknowledgements}

I would like to thank Anne Shiu and the anonymous referees for a number of helpful comments on the manuscript.

\appendix
\section{Methodological notes} 
\label{appmethod}
The processes of generating CRNs, and of searching numerically for oscillation, are described briefly. Further detail can be found in \cite{banajiCRNcount}.

{\bf Generating unlabelled CRNs.} Note first that two fully open CRNs are isomorphic if and only if they are isomorphic after removal from both of the flow reactions $0 \rightleftharpoons X_i$. All at most bimolecular unlabelled, fully open CRNs with $k$ species and $l$ non-flow reactions can be generated as follows: 
\begin{enumerate}[align=left,leftmargin=*]
\item All at most bimolecular complexes on $k$ species are listed: there are $n_C(k):={k+2 \choose 2}$ such complexes.
\item Irreversible reactions can be viewed as ordered pairs of distinct complexes: and consequently there are a total of $n_C(k)(n_C(k)-1)$ distinct irreversible reactions involving these complexes. 
\item Excluding the reactions $0 \rightarrow X_i$ and $X_i \rightarrow 0$ leaves $n_R(k):=n_C(k)(n_C(k)-1) - 2k$ distinct non-flow reactions from which to build the CRNs. 
\item All sets of $l$ distinct non-flow reactions are chosen and represented (in {\tt digraph6} format) as two-layer vertex-coloured digraphs, as described in the section {\em Isomorphism of edge-coloured graphs} of the NAUTY documentation at \url{http://users.cecs.anu.edu.au/~bdm/nauty/nug26.pdf}. There are
\[
{n_R(k) \choose l} = {{k+2 \choose 2}\left({k+2 \choose 2}-1\right) - 2k \choose l}
\]
of these CRNs. This number corresponds to the total number of $(k,l)$ CRNs with labelled species, but unlabelled reactions. 
\item The NAUTY program {\tt shortg} is used to canonically label and remove isomorphs from this list of CRNs, respecting the species-reaction bipartition. 
\end{enumerate}
Details of the enumeration methodology can be found in \cite{banajiCRNcount} with data at \url{https://reaction-networks.net/networks/}.

{\bf Generating $(k,l)$ CRNs which inherit oscillation.} Given lists of oscillatory $(k-1,l)$ CRNs and $(k,l-1)$ CRNs (with the reactions $0 \rightleftharpoons X_i$ removed), the following procedure generates $(k,l)$ CRNs which are oscillatory by inheritance:
\begin{enumerate}[align=left,leftmargin=*]
\item Each possible new non-flow reaction is added to each oscillatory $(k,l-1)$ CRN, giving a list of oscillatory $(k,l)$ CRNs.
\item The new species $X_{k}$ is added into the reactions of each oscillatory $(k-1,l)$ CRN in every possible way which preserves bimolecularity, giving a list of oscillatory $(k,l)$ CRNs.
\item The two lists obtained in the steps above are merged and {\tt shortg} is used to canonically label and remove isomorphs from the combined list. The open extensions of CRNs in the merged list are the $(k, l)$ CRNs which inherit stable oscillation, by Proposition~\ref{coropeninduced}. 
\end{enumerate}

{\bf Numerical simulations.} The following procedure was set up to search for oscillation in CRNs which neither inherited oscillation, nor were conjectured to forbid oscillation by Conjecture~\ref{conjJac}. For each such CRN:
\begin{enumerate}[align=left,leftmargin=*]
\item The differential equations were constructed from a combinatorial description of the CRN, along with the assumption of mass action kinetics or physical power-law kinetics. 
\item A minimum of 10,000 parameter-sets were chosen randomly using uniform distributions on each parameter. The parameters chosen were rate constants and initial conditions, and additionally exponents in the case of power-law kinetics. 

\item Simulations of the differential equations were carried out using {\tt RADAU5} software \cite{radau5}.
\item A script was written to analyse the outputs, searching for bounded but apparently nonconvergent trajectories. Where such behaviour was identified, plots of the trajectories were output and later
examined by eye to confirm that indeed oscillation had been found numerically. 
\end{enumerate} 
Several choices were necessarily somewhat arbitrary, particularly the number of simulations, the length of simulation, and the upper and lower limits on the magnitudes of parameters. It was also assumed throughout that what appeared in a plot as a periodic orbit was indeed an SPPO.

To explore the likelihood of finding oscillation in numerical simulations, from 100 to 100,000 simulations were carried out on each of the 444 $(3,3)$ CRNs known to admit SPPOs with mass action kinetics (see Table~\ref{tabdat}). The results, shown in Table~\ref{tabsimdat}, reflect the fact that oscillation often occurs only in small regions of parameter space, and so can be hard to find by brute-force approaches. In larger CRNs this problem becomes even more acute, highlighting the importance of theoretical approaches, including the inheritance results in this paper.
\begin{table}[h]
\centering
\begin{tabular}{|c|c|}
\cline{1-2}
parameter sets& CRNs found to admit SPPOs\\
\cline{1-2}
100 & 17\\
1,000 & 63\\
10,000 & 174\\
100,000 & 298\\
\cline{1-2}
\end{tabular}
\caption{\label{tabsimdat}Simulations were carried out on the 444 $(3,3)$ CRNs known to admit SPPOs with mass action kinetics to determine the effect of number of simulations on the likelihood of observing SPPOs in numerical experiment. For example, even when 100,000 simulations were carried out on each CRN, only 298 (67\%) of the CRNs were identified as oscillatory.}
\end{table}

\section{Proof of Proposition~\ref{prop21}}
\label{pf21}

\begin{myproof}{Proposition~\ref{prop21}}
The 14 non-isomorphic CRNs listed in the proposition are easily enumerated with NAUTY \cite{nauty} as described in \ref{appmethod} or even by eye. Under the assumption of positive general kinetics each of the CRNs $\mathcal{R}_{(i)}$ to $\mathcal{R}_{(xiv)}$ gives rise to an ODE system of the form
\begin{equation}
\label{gen2D}
\left.\begin{array}{rcl}
\dot x & = & a-g_1(x)+n_X f(x,y)\\
\dot y & = & c-g_2(y)+n_Y f(x,y)
\end{array}\right\}\,\,:=\,\,F(x,y)
\end{equation}
where $a$ and $c$ are positive constants; $g_1$ and $g_2$ are positive-valued $C^1$ functions on $\mathbb{R}_{>0}$ with positive derivative; and $f(x,y)$, the rate of the non-flow reaction, is a positive-valued $C^1$ function on $\mathbb{R}^2_{\gg 0}$ satisfying $\frac{\partial f}{\partial x} > 0$ (resp., $\frac{\partial f}{\partial y} > 0$) on $\mathbb{R}^2_{\gg 0}$ if $X$ (resp., $Y$) occurs on the left of the non-flow reaction. $n_X$ and $n_Y$ are the net production of $X$ and $Y$ respectively in the non-flow reaction. (\ref{gen2D}) defines a local flow $\phi$ on $\mathbb{R}^2_{\gg 0}$. 

(\ref{gen2D}) has Jacobian matrix
\[
J(x,y):=DF(x,y) = \left(\begin{array}{cc}-g_1'(x)+n_X f_x(x,y)&n_X f_y(x,y)\\n_Y f_x(x,y)&-g_2'(y)+n_Y f_y(x,y)\end{array}\right).
\]
on $\mathbb{R}^2_{\gg 0}$, and so
\[
\begin{array}{rcl}\mathrm{det}\,J(x,y) &=& g_1'(x)g_2'(y) -g_1'(x)n_Y f_y(x,y) -g_2'(y) n_X f_x(x,y)\,, \\\mathrm{Tr}\,J(x,y) &=& -g_1'(x)-g_2'(y) + n_X f_x(x,y) + n_Y f_y(x,y)\,.
\end{array}
\]
The assumption of positive general kinetics implies that $g_1'(x) > 0$ and $g_2'(y) > 0$ for all positive $x$ and $y$, and so $\mathrm{det}\,J(x,y)>0$ and $\mathrm{Tr}\,J(x,y)<0$ provided $n_X f_x(x,y)\leq 0$ and $n_Y f_y(x,y)\leq 0$. 

{\bf $\mathcal{R}_{(i)}$ to $\mathcal{R}_{(xii)}$.} In these cases, it is easily seen that the assumption of positive general kinetics ensures that $n_X f_x(x,y)\leq 0$ and $n_Y f_y(x,y)\leq 0$. Consequently, $\mathrm{Tr}\,J(x,y)<0$, and so each CRN forbids nontrivial periodic orbits in $\mathbb{R}^2_{\gg 0}$ by the Bendixson criterion (Theorem~4.1.1 of \cite{wiggins}). Notice that since $J(x,y)$ is everywhere Hurwitz stable on $\mathbb{R}^2_{\gg 0}$, all positive equilibria are locally asymptotically stable and, by the 2D Markus-Yamabe Theorem \cite{fessler}, existence of a positive equilibrium guarantees that it attracts all of $\mathbb{R}^2_{\gg 0}$. (To apply the Theorem as it is usually stated, we may first pass to logarithmic coordinates via $(x,y) \mapsto (\ln x,\ln y)$, a smooth diffeomorphism between $\mathbb{R}^2_{\gg 0}$ and $\mathbb{R}^2$.)

Under the assumption of mass action kinetics we can go further. We have $g_1(x) = bx$ and $g_2(y) = dy$ for some $b>0$ and $d>0$ and $f(x,y) = \gamma x^\alpha y^\beta$ where $\gamma>0$ and $\alpha$ and $\beta$ refer to the stoichiometries of $X$ and $Y$ on the left of the non-flow reaction. It is an easy exercise to prove that with mass action kinetics (and indeed more general assumptions), each of $\mathcal{R}_{(i)}$ to $\mathcal{R}_{(xii)}$ has a unique positive equilibrium which must, by the arguments above, attract all of $\mathbb{R}^2_{\gg 0}$. In fact, since each forward trajectory originating on $\partial \mathbb{R}^2_{\geq 0}$ immediately enters $\mathbb{R}^2_{\gg 0}$, the unique positive equilibrium attracts all of $\mathbb{R}^2_{\geq 0}$. As it is locally asymptotically stable, it is in fact globally asymptotically stable.

{\bf $\mathcal{R}_{(xiii)}$.} In this case, for positive general kinetics, (\ref{gen2D}) is a decoupled system of two ODEs which clearly forbids periodic orbits. With mass action kinetics, the system can be written
\begin{equation}
\label{gen2Dxiii}
\begin{array}{rcl}
\dot x & = & a-bx+\gamma x\\
\dot y & = & c-dy
\end{array}
\end{equation}
For $b>\gamma$, (\ref{gen2Dxiii}) has a unique, globally attracting, positive equilibrium at $(a/(b-\gamma), c/d)$, while for $b \leq \gamma$ it has only unbounded trajectories.

What remains is the only nontrivial case, namely $\mathcal{R}_{(xiv)}$:
\[
Y \rightleftharpoons 0 \rightleftharpoons X, \quad X+Y \rightarrow 2Y\,. 
\]

{\bf Existence of an SPPO for $\mathcal{R}_{(xiv)}$ for general kinetics and physical power-law kinetics.} The two cases are dealt with at once by proving the existence of an SPPO for reaction rates which belong simultaneously to both classes of kinetics.  Fixing mass-action kinetics for $Y \rightleftharpoons 0 \rightleftharpoons X$, but allowing physical power-law kinetics for the non-flow reaction, the system takes the form:
\begin{equation}
\label{xivPPLK}
\begin{array}{rcl}
\dot x & = & a-bx-\gamma x^\alpha y^\beta\\
\dot y & = & c-dy+\gamma x^\alpha y^\beta.
\end{array}
\end{equation}
where $\alpha>0$ and $\beta>0$. We now prove that (\ref{xivPPLK}) has an attracting periodic orbit for appropriate choices of the parameters. We choose
\[
a=\frac{3}{2},\,\, b=\frac{1}{2},\,\, c=\frac{1}{2}-k,\,\, d=\frac{3}{2}-k,\,\, \gamma = 1,\,\, \alpha=1\,\,\mbox{and}\,\,\beta=3.
\]
$k$ is a bifurcation parameter which affects two rate constants, and the choices ensure that the kinetics is polynomial and belongs both to the class of general kinetics and of physical power-law kinetics. With these parameters (\ref{xivPPLK}) becomes
\begin{equation}
\label{xivset}
\left.\begin{array}{rcl}
\dot x & = & \frac{3}{2}-\frac{x}{2}-x y^3\\
\dot y & = & \left(\frac{1}{2}-k\right)-\left(\frac{3}{2}-k\right)y+x y^3
\end{array}
\right\}=: F(x,y;k).
\end{equation}
Note that $F(1,1;k)=0$ for all $k$ and $DF(1,1;k)$ has eigenvalues
\[
\lambda_{\pm}(k):=\frac{k \pm \sqrt{k^2+6k-3}}{2}\,.
\]
We confirm that a Hopf bifurcation occurs at $k=0$. Let $\omega = \frac{\sqrt{3}}{2}$ so that $\lambda_{\pm}(0) = \pm i\omega$. A quick calculation gives that
\[
\left.\frac{\mathrm{d}}{\mathrm{d}k}(\mathrm{Re}\,\lambda_{\pm}(k))\right|_{k=0} = \frac{1}{2}>0\,,
\]
namely $\lambda_{\pm}(k)$ cross with nonzero speed from left half to right half plane as $k$ increases through $0$. To confirm that an attracting periodic orbit indeed exists for sufficiently small $k>0$ we bring (\ref{xivset}) with $k=0$ into a standard form by defining new variables $(u,v)$ via
\[
\left(\begin{array}{c}x\\y\end{array}\right) = \left(\begin{array}{rc}1&0\\-\frac{1}{2}&\frac{\sqrt{3}}{6}\end{array}\right)\left(\begin{array}{c}u\\v\end{array}\right)+\left(\begin{array}{c}1\\1\end{array}\right),
\]
to get
\begin{equation}
\label{xivseta}
\left(\begin{array}{c}
\dot u\\
\dot v \end{array}\right) = \left(\begin{array}{cc}0&-\omega\\\omega&0\end{array}\right)\left(\begin{array}{c}u\\v\end{array}\right) + \left(\begin{array}{c}f^1(u,v)\\f^2(u,v)\end{array}\right).
\end{equation}
Here (\ref{xivseta}) has an equilibrium at $(0,0)$ and $f^1(u,v)$, $f^2(u,v)$ consist of the nonlinear terms in $u,v$, i.e., $f^1(0,0)=f^2(0,0)=f^1_u(0,0)=f^1_v(0,0)=f^2_u(0,0)=f^2_v(0,0)=0$. We can calculate the relevant combination of higher partial derivatives of $f^1$ and $f^2$ at $(0,0)$ to get:
\[
\frac{1}{16}(f^1_{uuu}+f^1_{uvv}+f^2_{uuv}+f^2_{vvv})+\frac{1}{16\omega}(f^1_{uv}(f^1_{uu}+f^1_{vv})-f^2_{uv}(f^2_{uu}+f^2_{vv}) - f^1_{uu}f^2_{uu}+f^1_{vv}f^2_{vv})=-\frac{1}{8}<0.
\]
This computation, carried out with the help of MAXIMA \cite{maxima}, implies that as $k$ increases through zero we have a supercritical Hopf bifurcation and the creation of an asymptotically orbitally stable periodic orbit near to the origin in $u\mhyphen v$ space (see Section~2 of Chapter~20 in \cite{wiggins}). In practice, stable periodic orbits are easily found in numerical simulations of (\ref{xivPPLK}) for many choices of parameters. 

{\bf Nonexistence of a periodic orbit for $\mathcal{R}_{(xiv)}$ with mass action kinetics.} We now show that with mass action kinetics, $\mathcal{R}_{(xiv)}$ admits no periodic orbits, and in fact has a unique positive equilibrium which is globally asymptotically stable. Specialising to mass action kinetics, (\ref{xivPPLK}) becomes
\begin{eqnarray*}
\dot x & = & a-bx-\gamma xy\\
\dot y & = & c-dy+\gamma xy.
\end{eqnarray*}
Rescaling $x \mapsto \gamma x$, $y \mapsto \gamma y$, $a\mapsto \gamma a$ and $c \mapsto \gamma c$ allows us to eliminate $\gamma$ from the system to get:
\begin{equation}
\label{NOPO}
\left.\begin{array}{rcl}
\dot x & = & a-bx-xy\\
\dot y & = & c-dy+xy.
\end{array}\right\}\,\,:=\,\,F(x,y)\,.
\end{equation}
(\ref{NOPO}) defines a local semiflow $\phi$ on $\mathbb{R}^2_{\geq 0}$. Let $G(x,y):= \frac{1}{y}F(x,y)$. Then
\[
\mathrm{Tr}DG(x,y) = \frac{1}{y}(-b-y-d+x) - \frac{1}{y^2}(c-dy+xy) = \frac{1}{y}\left(-b-y -\frac{c}{y}\right) < 0\,.
\]
As this holds everywhere on $\mathbb{R}^2_{\gg 0}$, (\ref{NOPO}) has no nontrivial periodic orbits in $\mathbb{R}^2_{\gg 0}$ by the Bendixson-Dulac criterion (Theorem~4.1.2 of \cite{wiggins}). 

Next, we show that (\ref{NOPO}) has a unique linearly stable equilibrium which attracts all of $\mathbb{R}^2_{\geq 0}$. As each point of $\partial \mathbb{R}^2_{\geq 0}$ is a start point of $\phi$, $\phi$ has no $\omega$-limit points on $\partial \mathbb{R}^2_{\geq 0}$ (Theorem~5.6 in \cite{bhatiahajek} for example). Note that the triangle $T$ defined by
\[
x \geq 0, \quad y \geq 0, \quad x+y \leq \frac{(b+d)(a+c)}{bd}\,
\]
is a global attractor for (\ref{NOPO}) as
\[
\frac{\mathrm{d}}{\mathrm{d} t}(x+y) = a+c-bx-dy < a+c -\frac{bd}{b+d}(x+y) \leq 0 \quad \mbox{provided} \quad x+y \geq \frac{(b+d)(a+c)}{bd}.
\]
Thus all orbits of $\phi$ are bounded and $\phi$ is in fact a semiflow on $\mathbb{R}^2_{\geq 0}$. Moreover each point of $\partial T$ is a start point of $\left.\phi\right|_{T}$, and so $\left.\phi\right|_{T}$ (and hence $\phi$) has no $\omega$-limit points on $\partial T$. 

There is a unique positive solution $(x_0,y_0)$ to $F(x,y)=0$ given by
\[
\begin{array}{rcl}x_0 &=& \displaystyle{\frac{a+c+bd-\theta}{2b}}\\
y_0 &=& \displaystyle{\frac{a+c-bd+\theta}{2d}}
\end{array} \quad \mbox{where} \quad \theta := \sqrt{(a+c+bd)^2 - 4abd} = \sqrt{(a+c-bd)^2 + 4bcd}\,.
\]
As $\theta > a+c-bd$, we see that $x_0 < d$. The Jacobian matrix of the system is
\[
J = \left(\begin{array}{cc}-b-y&-x\\y&-d+x\end{array}\right).
\]
At $(x_0,y_0)$, as $x_0 < d$,
\[
\mathrm{det}J(x_0,y_0) = b(d-x_0) + dy_0 > 0\quad\mbox{and}\quad \mathrm{Tr}J(x_0,y_0) = -b-y_0-(d-x_0) < 0\,.
\]
and so $(x_0, y_0)$ is linearly stable. Certainly, $(x_0,y_0)$ has no homoclinic orbits. Since all trajectories of (\ref{NOPO}) enter $\mathrm{int}\,T$ and $\phi$ has no nontrivial periodic or homoclinic orbits in $\mathrm{int}\,T$, by the Poincar\'e-Bendixson theorem (Theorem~9.0.6 in \cite{wiggins} for example), $(x_0,y_0)$ is a global attractor of (\ref{NOPO}) and hence is globally asymptotically stable.
\hfill
\end{myproof}

\bibliographystyle{unsrt}

\begin{thebibliography}{10}

\bibitem{MurrayMathBio}
J.~D. Murray.
\newblock {\em Mathematical Biology.}
\newblock Springer, Berlin, 1993.

\bibitem{mathphys}
J.~Keener and J.~Sneyd.
\newblock {\em Mathematical Physiology}, volume~8 of {\em Interdisciplinary
  Applied Mathematics}.
\newblock Springer, New York, 1998.

\bibitem{novaktyson}
B.~Nov\'ak and J.J. Tyson.
\newblock Design principles of biochemical oscillators.
\newblock {\em Nature Rev. Mol. Cell. Biol.}, 9:981--991, 2008.

\bibitem{joshishiu}
B.~Joshi and A.~Shiu.
\newblock Atoms of multistationarity in chemical reaction networks.
\newblock {\em J. Math. Chem.}, 51(1):153--178, 2013.

\bibitem{feliuwiufInterface2013}
E.~Feliu and C.~Wiuf.
\newblock Simplifying biochemical models with intermediate species.
\newblock {\em J. Roy. Soc. Interface}, 10:20130484, 2013.

\bibitem{Joshi.2013aa}
B.~Joshi.
\newblock Complete characterization by multistationarity of fully open networks
  with one non-flow reaction.
\newblock {\em Applied Mathematics and Computation}, 219(12):6931--6945, 2013.

\bibitem{JoshiShiu2016}
B.~Joshi and A.~Shiu.
\newblock Which small reaction networks are multistationary?
\newblock {\em SIAM J. Appl. Dyn. Syst.}, 16(2):802--833, 2017.

\bibitem{banajipanteaMPNE}
M.~Banaji and C.~Pantea.
\newblock The inheritance of nondegenerate multistationarity in chemical
  reaction networks.
\newblock preprint at {\tt https://arxiv.org/abs/1608.08400}.

\bibitem{horn72}
F.~Horn.
\newblock Necessary and sufficient conditions for complex balancing in chemical
  kinetics.
\newblock {\em Arch. Ration. Mech. Anal.}, 49:172--186, 1972.

\bibitem{hornjackson}
F.~Horn and R.~Jackson.
\newblock General mass action kinetics.
\newblock {\em Arch. Ration. Mech. Anal.}, 47(2):81--116, 1972.

\bibitem{feinberg0}
M.~Feinberg.
\newblock Complex balancing in general kinetic systems.
\newblock {\em Arch. Ration. Mech. Anal.}, 49(3):187--194, 1972.

\bibitem{feinberg}
M.~Feinberg.
\newblock Chemical reaction network structure and the stability of complex
  isothermal reactors - {I}. {T}he deficiency zero and deficiency one theorems.
\newblock {\em Chem. Eng. Sci.}, 42(10):2229--2268, 1987.

\bibitem{banajidynsys}
M.~Banaji.
\newblock Monotonicity in chemical reaction systems.
\newblock {\em Dyn. Syst.}, 24(1):1--30, 2009.

\bibitem{angelileenheersontag}
D.~Angeli, P.~{De Leenheer}, and E.~D. Sontag.
\newblock Graph-theoretic characterizations of monotonicity of chemical
  reaction networks in reaction coordinates.
\newblock {\em J. Math. Biol.}, 61(4):581--616, 2010.

\bibitem{donnellbanaji}
P.~Donnell and M.~Banaji.
\newblock Local and global stability of equilibria for a class of chemical
  reaction networks.
\newblock {\em SIAM J. Appl. Dyn. Syst.}, 12(2):899--920, 2013.

\bibitem{banajimierczynski}
M.~Banaji and J.~Mierczy\'nski.
\newblock Global convergence in systems of differential equations arising from
  chemical reaction networks.
\newblock {\em J. Diff. Eq.}, 254(3):1359--1374, 2013.

\bibitem{abphopf}
D.~Angeli, M.~Banaji, and C.~Pantea.
\newblock Combinatorial approaches to {H}opf bifurcations in systems of
  interacting elements.
\newblock {\em Commun Math Sci}, 12:1101--1133, 2014.

\bibitem{angelihirschsontag}
D.~Angeli, M.~W. Hirsch, and E.~Sontag.
\newblock Attractors in coherent systems of differential equations.
\newblock {\em J. Diff. Eq.}, 246:3058--3076, 2009.

\bibitem{li_muldowney_1993}
M.~Y. Li and J.~S. Muldowney.
\newblock On {Bendixson's} criterion.
\newblock {\em J. Diff. Eq.}, 106:27--39, 1993.

\bibitem{li_muldowney_1996}
M.~Y. Li and J.~S. Muldowney.
\newblock A geometric approach to global-stability problems.
\newblock {\em SIAM J. Math. Anal.}, 27(4):1070--1083, 1996.

\bibitem{li_muldowney_2000}
M.~Y. Li and J.~S. Muldowney.
\newblock Dynamics of differential equations on invariant manifolds.
\newblock {\em J. Diff. Eq.}, 168:295--320, 2000.

\bibitem{dicera}
E.~{Di Cera}, P.E. {Phillipson}, and J.~Wyman.
\newblock Limit-cycle oscillations and chaos in reaction networks subject to
  conservation of mass.
\newblock {\em Proc. Natl. Acad. Sci. USA}, 86:142--146, 1989.

\bibitem{WolfOsci}
J.~Wolf, {H-Y.} Sohn, R.~Heinrich, and H.~Kuriyama.
\newblock Mathematical analysis of a mechanism for autonomous metabolic
  oscillations in continuous culture of {Saccharomyces cerevisiae}.
\newblock {\em FEBS Letters}, 499(3):230--234, 2001.

\bibitem{Kholodenko.2000aa}
B.~N. Kholodenko.
\newblock Negative feedback and ultrasensitivity can bring about oscillations
  in the mitogen-activated protein kinase cascades.
\newblock {\em Eur. J. Biochem.}, 267(6):1583--1588, 2000.

\bibitem{Qiao.2007aa}
L.~Qiao, R.~B. Nachbar, I.~G. Kevrekidis, and S.~Y. Shvartsman.
\newblock {Bistability and oscillations in the Huang-Ferrell model of MAPK
  signaling}.
\newblock {\em PLoS Comput. Biol}, pages 1819--1826, 2007.

\bibitem{eiswirth91}
M.~Eiswirth, A.~Freund, and J.~Ross.
\newblock Mechanistic classification of chemical oscillators and the role of
  species.
\newblock In {\em Adv. Chem. Phys. {Vol 80}}, chapter~2. John Wiley \& sons,
  1991.

\bibitem{eiswirth96}
M.~Eiswirth, J.~B\"urger, P.~Strasser, and G.~Ertl.
\newblock {Oscillating Langmuir-Hinshelwood mechanisms}.
\newblock {\em J. Phys. Chem.}, 100(49):19118--19123, 1996.

\bibitem{gatermann}
K.~Gatermann, M.~Eiswirth, and A.~E. Sensse.
\newblock Toric ideals and graph theory to analyze {H}opf bifurcations in mass
  action systems.
\newblock {\em J. Symbolic Comput.}, 40:1361--1382, 2005.

\bibitem{errami2015}
H.~Errami, M.~Eiswirth, D.~Grigoriev, W.~M. Seiler, T.~Sturm, and A.~Weber.
\newblock {Detection of Hopf bifurcations in chemical reaction networks using
  convex coordinates}.
\newblock {\em J. Comput. Phys.}, 291:279--302, 2015.

\bibitem{minchevaroussel}
M.~Mincheva and M.~R. Roussel.
\newblock Graph-theoretic methods for the analysis of chemical and biochemical
  networks, {I.} {M}ultistability and oscillations in ordinary differential
  equation models.
\newblock {\em J. Math. Biol.}, 55:61--86, 2007.

\bibitem{domijan}
M.~Domijan and M.~Kirkilionis.
\newblock Bistability and oscillations in chemical reaction networks.
\newblock {\em J. Math. Biol.}, 59:467--501, 2009.

\bibitem{ConradiShiuPTM}
C.~Conradi and A.~Shiu.
\newblock Dynamics of post-translational modification systems: recent progress
  and future directions.
\newblock preprint at {\tt https://arxiv.org/abs/1705.10913}.

\bibitem{hadac}
O.~Hada\v{c}, F.~Muzika, V.~Nevoral, M.~P\v{r}ibyl, and I.~Schreiber.
\newblock {Minimal oscillating subnetwork in the Huang-Ferrell model of the
  MAPK cascade}.
\newblock {\em PLoS One}, 12(6):e0178457, 2017.

\bibitem{HaleOsci}
Jack~K. Hale.
\newblock {\em Oscillations in Nonlinear Systems}.
\newblock Dover, New York, 1963.

\bibitem{wiggins}
S.~Wiggins.
\newblock {\em Introduction to Applied Nonlinear Dynamics and Chaos}.
\newblock Springer, 2003.

\bibitem{Fenichel79}
N.~Fenichel.
\newblock Geometric singular perturbation theory for ordinary differential
  equations.
\newblock {\em J. Differ. Equations}, 31:53--98, 1979.

\bibitem{banajipantea}
M.~Banaji and C.~Pantea.
\newblock Some results on injectivity and multistationarity in chemical
  reaction networks.
\newblock {\em SIAM J. Appl. Dyn. Syst.}, 15(2):807--869, 2016.

\bibitem{horn}
F.~Horn.
\newblock On a connexion between stability and graphs in chemical kinetics,
  {I}. {S}tability and the reaction diagram, {II}. {S}tability and the complex
  graph.
\newblock {\em Proc. Roy. Soc. Lond. A}, 334:299--330, 1973.

\bibitem{angelipetrinet}
D.~Angeli, P.~{De Leenheer}, and E.~D. Sontag.
\newblock A {P}etri net approach to the study of persistence in chemical
  reaction networks.
\newblock {\em Math. Biosci.}, 210:598--618, 2007.

\bibitem{banajiCRNcount}
M.~Banaji.
\newblock {Counting chemical reaction networks with NAUTY}.
\newblock \url{https://arxiv.org/abs/1705.10820}.

\bibitem{bhatiahajek}
N.~P. Bhatia and O.~Hajek.
\newblock {\em Local semi-dynamical systems}.
\newblock Springer-Verlag, 1969.

\bibitem{kuznetsov}
Y.~A. Kuznetsov.
\newblock {\em Elements of applied bifurcation theory}.
\newblock Springer, New York, 1998.

\bibitem{nauty}
B.~D. McKay and A.~Piperno.
\newblock Practical graph isomorphism {II}.
\newblock {\em J. Symbolic Computation}, 60:94--112, 2013.

\bibitem{pojman}
I.~R. Epstein and J.~A. Pojman, editors.
\newblock {\em An Introduction to Nonlinear Chemical Dynamics: Oscillations,
  Waves, Patterns, and Chaos}.
\newblock Oxford University Press, New York, 1998.

\bibitem{radau5}
E.~Hairer and G.~Wanner.
\newblock {RADAU5} code for integrating differential algebraic equations.
\newblock Available online at
  \url{http://www.unige.ch/~hairer/prog/stiff/radau5.f}.

\bibitem{fessler}
R.~Fe{\ss}ler.
\newblock A proof of the two-dimensional {M}arkus-{Y}amabe stability
  conjecture.
\newblock {\em Ann. Polon. Math.}, 62:45--75, 1995.

\bibitem{maxima}
{MAXIMA}: A computer algebra system.
\newblock Available at http://maxima.sourceforge.net/.

\end{thebibliography}

\end{document}